\documentclass[a4paper,fleqn]{cas-dc}

\usepackage[numbers]{natbib}

\usepackage{amsmath,amssymb,amsfonts}
\usepackage{algorithm,algorithmic}
\usepackage{booktabs}

\newproof{proof}{Proof}

\newif\ifcomments\commentsfalse

\begin{document}
\let\WriteBookmarks\relax
\def\floatpagepagefraction{1}
\def\textpagefraction{.001}

\shorttitle{}    

\shortauthors{L.~L\"uken, S.~Lucia}  

\title [mode = title]{A Scalable Learned Iterative Solver for Nonlinear MPC via Temporal Decomposition}  

\author[1]{Lukas~L\"uken}[orcid=0009-0004-5599-8282]
\cormark[1]
\credit{Writing - original draft, Investigation, Visualization, Software, Methodology, Formal analysis, Conceptualization}
\ead{lukas.lueken@tu-dortmund.de}
\cortext[1]{Corresponding author}

\author[1]{Sergio~Lucia}[orcid=0000-0002-3347-5593]
\credit{Writing - review and editing, Supervision, Conceptualization}
\ead{sergio.lucia@tu-dortmund.de}

\affiliation[1]{organization={Chair of Process Automation Systems, TU Dortmund University},
            addressline={Emil-Figge Straße 70}, 
            city={Dortmund},
            postcode={44227}, 
            state={North Rhine-Westphalia},
            country={Germany}}

\begin{abstract}
Nonlinear model predictive control (NMPC) can be applied to large nonlinear systems with complex objectives while enforcing constraints, but the cost of solving the underlying optimization problems limits real-time application for systems with fast dynamics or many states.

Neural networks can reduce this online cost by approximating the control law, but may violate constraints and typically require many precomputed optimizer solutions for training.
The learning-based iterative solver for constrained optimization (LISCO) instead learns to predict primal-dual steps toward satisfaction of the Karush-Kuhn-Tucker (KKT) conditions through self-supervised training without precomputed solutions, but its dense, monolithic architecture grows with the prediction horizon.

We remove this horizon dependence by exploiting the temporal structure of the problem in the solver architecture, which makes the network size independent of the horizon while leaving the training loss and the solver iteration unchanged.
A closed-loop resampling strategy further aligns the training distribution with closed-loop operation at no additional sampling cost.

On an economic NMPC benchmark of a CSTR cascade, the proposed solver scales to larger problems than monolithic LISCO and achieves median speedups of up to a factor of four over IPOPT.
In closed loop, it achieves comparable control performance to IPOPT with no observed constraint violations.

\end{abstract}


\maketitle

\section{Introduction}\label{sec:introduction}
Nonlinear model predictive control (NMPC) is highly relevant in process control, as it optimizes performance directly while enforcing input and state constraints \cite{ellisTutorialReviewEconomic2014,rawlingsFundamentalsEconomicModel2012a}.
Its online application requires solving a nonlinear optimization problem at every sampling instant.
The computational cost of this problem grows with the prediction horizon, the model size and the number of constraints \cite{rawlingsModelPredictiveControl2017}.
Real-time application is therefore limited for systems with fast dynamics, for large systems with many states and inputs, and for robust NMPC formulations \cite{luciaMultistageNonlinearModel2013}.

The online cost of NMPC has been addressed with efficient formulations such as the real-time iteration scheme \cite{diehlRealTimeIterationScheme2005} and fast solvers which exploit the structure of the NMPC problem \cite{frisonHPIPMHighperformanceQuadratic2020,verschuerenAcadosModularOpensource2022,zanelli_forces_2020,ferreauEmbeddedOptimizationMethods2017}.
These methods lower the online cost substantially, but still require solving structured linear systems at each iteration, with a computational cost that grows with the problem size.
For small linear problems, explicit MPC precomputes the control law offline \cite{bemporadExplicitLinearQuadratic2002}, which does not extend to large nonlinear problems.
This motivates learning-based methods, which are fast to evaluate and easy to deploy on control hardware.

Existing learning-based methods can be broadly categorized into two areas, the direct approximation of the solution and the learning of the solver itself.
In the first area, a neural network directly predicts the optimization solution or the control input.
Approximate MPC follows this approach by approximating the MPC control law, which maps the current state to the first control input of the optimal input sequence \cite{parisiniRecedinghorizonRegulatorNonlinear1995,kargEfficientRepresentationApproximation2020,lukenSobolevTrainingDataefficient2023,chenApproximatingExplicitModel2018}.
Evaluating the network is often orders of magnitude faster than solving the problem online, but approximation errors can cause constraint violations and suboptimal closed-loop behavior \cite{hertneckLearningApproximateModel2018,kargProbabilisticPerformanceValidation2021}, and supervised training typically requires sampling large numbers of optimizer solutions, which is computationally expensive.
For constrained optimization more generally, the framework DC3 \cite{dontiDC3LearningMethod2021} adds a completion and correction step that reduces constraint violations, but incurs a computational overhead and leaves a residual optimality gap.
Primal-dual learning and related self-supervised optimization proxies \cite{parkSelfsupervisedPrimaldualLearning2023,fiorettoLagrangianDualityConstrained2020,kotaryEndtoEndConstrainedOptimization2021,frommeUnsupervisedClosedloopPrimaldual2025} remove the need for sampled optimizer solutions.
However, their direct predictions can still be infeasible or suboptimal.
A further line of work learns an initial guess that warm-starts a classical solver \cite{klaucoMachineLearningbasedWarm2019,vaupelAcceleratingNonlinearModel2020,sambharyaEndtoendLearningWarmstart2023}, which accelerates the solution but still requires an online solver.

In the second area, an update map or algorithm is learned and applied iteratively, so that the residual can be recomputed and monitored online.
In the learning-to-optimize (L2O) literature, the iteration map itself is learned, mostly for unconstrained problems such as neural network training \cite{chenLearningOptimizePrimer2022,liLearningOptimize2016}.
Related approaches tune or optimize the hyperparameters of a given algorithm \cite{mitsosOptimalDeterministicAlgorithm2018,sambharya2026learning} or unroll it and fit the parameters of the unrolled scheme \cite{shlezingerModelBasedDeepLearning2023,liPDHGUnrolledLearningtoOptimizeMethod2024}.
In these cases, learning builds on a prescribed optimization algorithm to accelerate the solution process.
In our previous work, we introduced the learning-based iterative solver for constrained optimization (LISCO), tailored to constrained nonlinear programs (NLP) and applied to NMPC \cite{lukenSelfsupervisedLearningIterative2026,lukenLearningIterativeSolvers2024}.
It learns a primal-dual update map for a modified Karush-Kuhn-Tucker (KKT) residual, is trained self-supervised without sampled optimizer solutions, and produces full primal-dual iterates that support a residual-based stopping criterion.
Each iteration is a single network evaluation and therefore requires no linear system solve.
This approach has since been applied to mixed-integer and robust NMPC \cite{adamekEnablingRobustMixedinteger2026}.

However, LISCO remains monolithic, as one dense network maps the full primal-dual residual to the full step.
The input and output layers of such a network grow with the primal-dual dimension and therefore with the horizon, so the network size and the training effort have to grow with the problem \cite{lukenSelfsupervisedLearningIterative2026}.
Consistent with this, in our previous work we observed decreasing speedups over IPOPT \cite{wachterImplementationInteriorpointFilter2006} as the horizon length increased, with larger networks and more training epochs required to retain convergence \cite{lukenSelfsupervisedLearningIterative2026}.
Classical solvers such as IPOPT exploit the block-banded KKT structure induced by the dynamics and stage-wise costs and constraints using sparse linear algebra, e.g. MA27 \cite{duff1983multifrontal}, whereas the learned solver has so far not made use of this temporal structure.
In this work, we close this gap by incorporating the temporal structure of the NMPC problem into the learned update map.

Classical decomposition methods, such as the alternating direction method of multipliers \cite{boydDistributedOptimizationStatistical2011} and overlapping Schwarz schemes \cite{naSuperconvergenceOnlineOptimization2023,shinExponentialDecaySensitivity2022,bartaliSchwarzDecompositionParallel2024}, also exploit the temporal locality of optimal control problems, but rely on local subproblem solves, multiplier updates or explicit coordination at solve time.
The proposed solver handles this coupling through overlapping residual windows and overlap averaging within the learned update map, so the solver iteration and the self-supervised training loss remain unchanged.

The contributions of this paper are as follows.
\begin{enumerate}
    \item We introduce a learned primal-dual solver for NMPC based on temporal decomposition.
    A shared neural network computes local updates in parallel on overlapping segments of the prediction horizon.
    Combining these local updates yields a full primal-dual step, while sharing the network across segments makes its size independent of the prediction horizon.
    \item We extend the self-supervised training by a closed-loop resampling strategy, in which the initial state of a converged training sample is replaced by the predicted next state instead of being resampled uniformly.
    Since this state is already contained in the converged primal-dual iterate, the training distribution is aligned with closed-loop operation at no additional sampling cost.
\end{enumerate}

The remainder of this paper is structured as follows.
We introduce the NMPC problem, its formulation as a parametric NLP with the modified KKT residual, the monolithic LISCO solver and the block-banded structure of the NMPC KKT system in Section~\ref{sec:background}.
Next, we present the proposed architecture and its training in Section~\ref{sec:method}, describe the CSTR-cascade benchmark and the computational studies in Section~\ref{sec:case_study}, report the results in Section~\ref{sec:results}, and conclude in Section~\ref{sec:conclusion}.

\section{Background}\label{sec:background}

\subsection{Nonlinear Model Predictive Control}
\label{sec:bg_nmpc}

Consider a discrete-time system $\mathbf{x}_{i+1} = f(\mathbf{x}_i,\mathbf{u}_i)$ with initial condition $\mathbf{x}_0 = \mathbf{x}_{\mathrm{init}}$, states $\mathbf{x}_i \in \mathbb{R}^{n_x}$, inputs $\mathbf{u}_i \in \mathbb{R}^{n_u}$, dynamics $f$ and time index $i \in \{0,\ldots,N-1 ~|~ N \in \mathbb{N},\, N \geq 1\}$.
The input sequence is determined by solving the finite-horizon NMPC problem \eqref{eq:nmpc_problem}, with stage cost $\ell$, terminal cost $V_f$, the discretized dynamics as equality constraints, and inequality constraints $\mathbf{g}_{\mathrm{MPC}}$ \cite{rawlingsModelPredictiveControl2017}:
\begin{subequations}\label{eq:nmpc_problem}
    \begin{align}
        \min_{\mathbf{x}_{[0:N]},\, \mathbf{u}_{[0:N-1]}} \quad
            & V_f\left(\mathbf{x}_N\right) + \sum_{i=0}^{N-1} \ell\left(\mathbf{x}_i, \mathbf{u}_i\right)\\
        \mathrm{s.t.} \quad
            & \mathbf{x}_{i+1} = f\left(\mathbf{x}_i, \mathbf{u}_i\right), \quad \mathbf{x}_0 = \mathbf{x}_{\mathrm{init}},\\
            & \mathbf{g}_{\mathrm{MPC}}\left(\mathbf{x}_{i+1}, \mathbf{u}_i\right) \leq \mathbf{0},\\
            & i \in \{0,\ldots,N-1 ~|~ N \in \mathbb{N},\, N \geq 1\}.
    \end{align}
\end{subequations}
Following the receding-horizon principle, the problem is solved at the current state $\mathbf{x}_0$, the first input $\mathbf{u}_0$ is applied, and the procedure is repeated at the next sampling instant, which defines the implicit control law: 
\begin{equation}\label{eq:nmpc_control_law}
    \mathbf{u}_0 = \Pi_{\mathrm{MPC}}(\mathbf{x}_{\mathrm{init}}).
\end{equation}

A direct alternative to solving \eqref{eq:nmpc_problem} online is to approximate the control law itself.
Approximate MPC fits a neural network $\Pi_{\mathrm{approxMPC}}$ with parameters $\theta$ to the map from the initial state to the first control action, using precomputed optimizer solutions as training targets, so that evaluating the network online replaces solving the optimization problem \cite{kargEfficientRepresentationApproximation2020,hertneckLearningApproximateModel2018}:
\begin{equation}\label{eq:approx_control_law}
    \hat{\mathbf{u}}_0 = \Pi_{\mathrm{approxMPC}}(\mathbf{x}_{\mathrm{init}};\theta).
\end{equation}
It is often orders of magnitude faster than solving \eqref{eq:nmpc_problem} online, but leads to approximation errors that may result in constraint violations or loss in performance.

\subsection{Parametric Nonlinear Optimization and Optimality Conditions}
\label{sec:bg_nlp}

The NMPC problem \eqref{eq:nmpc_problem} can be rewritten as a parametric nonlinear program \eqref{eq:nlp} with objective $q$, decision variables $\mathbf{w} \in \mathbb{R}^{n_w}$, parameters $\mathbf{p} \in \mathbb{R}^{n_p}$, inequality constraints $\mathbf{g} \in \mathbb{R}^{n_g}$ and equality constraints $\mathbf{h} \in \mathbb{R}^{n_h}$, all assumed smooth and possibly nonconvex:
\begin{subequations}\label{eq:nlp}
    \begin{align}
        \min_{\mathbf{w}} \quad & q\left(\mathbf{w},\mathbf{p}\right)\\
        \mathrm{s.t.} \quad & \mathbf{g}\left(\mathbf{w},\mathbf{p}\right) \leq \mathbf{0},\\
        & \mathbf{h}\left(\mathbf{w},\mathbf{p}\right) = \mathbf{0}.
    \end{align}
\end{subequations}
Collecting the primal variables and the equality and inequality multipliers into the primal-dual vector $\mathbf{z} = (\mathbf{w},\boldsymbol{\nu},\boldsymbol{\lambda}) \in \mathbb{R}^{n_z}$ and writing the Lagrangian $\mathcal{L}(\mathbf{z},\mathbf{p}) = q(\mathbf{w},\mathbf{p}) + \boldsymbol{\nu}^{\top}\mathbf{h}(\mathbf{w},\mathbf{p}) + \boldsymbol{\lambda}^{\top}\mathbf{g}(\mathbf{w},\mathbf{p})$, the first-order optimality conditions are the KKT conditions \eqref{eq:kkt_conditions}:
\begin{subequations}\label{eq:kkt_conditions}
    \begin{align}
        \nabla_{\mathbf{w}}\mathcal{L}\left(\mathbf{z}^{*},\mathbf{p}\right) &= \mathbf{0}, \label{subeq:stationarity_condition}\\
        \mathbf{h}\left(\mathbf{w}^{*},\mathbf{p}\right) &= \mathbf{0}, \label{subeq:equality_constraint}\\
        \mathbf{g}\left(\mathbf{w}^{*},\mathbf{p}\right) &\leq \mathbf{0}, \label{subeq:inequality_constraint}\\
        \boldsymbol{\lambda}^{*} &\geq \mathbf{0}, \label{subeq:dual_feasibility}\\
        \boldsymbol{\lambda}^{*} \odot \mathbf{g}\left(\mathbf{w}^{*},\mathbf{p}\right) &= \mathbf{0}. \label{subeq:complementary_slackness}
    \end{align}
\end{subequations}
To solve the parametric NLP \eqref{eq:nlp} for a given parameter $\mathbf{p}$, iterative solvers are used to drive the KKT residual to a given tolerance \cite{nocedalNumericalOptimization2006}.
These iterative solvers are generally described by the following iterative update rule:
\begin{equation}\label{eq:iterative_solver}
    \mathbf{z}_{k+1} = \mathbf{z}_k + \alpha_k \Delta\mathbf{z}_k ,
\end{equation}
with iterate $\mathbf{z}_k$, step $\Delta\mathbf{z}_k$ and step size $\alpha_k$.

The NMPC problem \eqref{eq:nmpc_problem} is an instance of \eqref{eq:nlp} under the identification $\mathbf{w} = (\mathbf{u}_0,\mathbf{x}_1,\ldots,\mathbf{u}_{N-1},\mathbf{x}_N)$ and $\mathbf{p} = \mathbf{x}_{\mathrm{init}}$, where the state and input trajectories are the decision variables, the initial state is the only parameter, the discretized dynamics are the equality constraints, and the stage constraints are the inequality constraints.

\subsection{Modified KKT Residual}
\label{sec:bg_residual}

The KKT conditions \eqref{eq:kkt_conditions} contain inequalities~\eqref{subeq:inequality_constraint}, \eqref{subeq:dual_feasibility} and the complementarity condition~\eqref{subeq:complementary_slackness}, which hinder the direct application of root-finding methods for nonlinear systems of equations that drive the residual to zero.
The penalized Fischer-Burmeister function \eqref{eq:fb} replaces the inequality and complementarity conditions by a single equation per inequality constraint \cite{fischerSpecialNewtontypeOptimization1992,chenPenalizedFischerBurmeisterNCPfunction2000,liao-mcphersonFBstabProximallyStabilized2020}:
\begin{align}\label{eq:fb}
    \phi_j(\lambda_j, g_j) =\;& \rho_{\mathrm{FB}}\left(\lambda_j - g_j - \sqrt{\lambda_j^2 + g_j^2 + \epsilon^2}\right) \notag\\
    &+ (1-\rho_{\mathrm{FB}})\,\lambda_j^{+}\, g_j^{-},
    \qquad j = 1,\ldots,n_g ,
\end{align}
where $\lambda_j^{+} = \max\{\lambda_j,0\}$ and $g_j^{-} = \max\{-g_j,0\}$ and $\rho_{\mathrm{FB}} \in (0,1]$ is the penalization parameter and $\epsilon > 0$ the smoothing parameter to improve conditioning.
For $\epsilon > 0$, the Fischer-Burmeister term is smooth, while the penalty term remains nonsmooth for $\rho_{\mathrm{FB}} < 1$.
For any $\epsilon \geq 0$, satisfaction of $\phi_j(\lambda_j, g_j) = 0$ ensures that $g_j \leq 0$ and $\lambda_j \geq 0$ and $\epsilon$ acts as a slack on the complementarity condition such that $\lvert\lambda_j g_j\rvert = \epsilon^2/2$ for $\rho_{\mathrm{FB}}=1$ and $\lvert\lambda_j g_j\rvert \leq \epsilon^2/2$ for $\rho_{\mathrm{FB}} < 1$.

Replacing the inequality, dual-feasibility and complementarity conditions with $\boldsymbol{\phi} = \mathbf{0}$ yields the modified KKT residual \eqref{eq:r_phi}, whose zeros recover the original KKT conditions in the limit $\epsilon \rightarrow 0$:
\begin{equation}\label{eq:r_phi}
    \mathbf{r}_{\phi}(\mathbf{z},\mathbf{p}) =
    \begin{bmatrix}
        \nabla_{\mathbf{w}}\mathcal{L}\left(\mathbf{z},\mathbf{p}\right)\\
        \mathbf{h}\left(\mathbf{w},\mathbf{p}\right)\\
        \boldsymbol{\phi}\left(\boldsymbol{\lambda},\mathbf{g}\left(\mathbf{w},\mathbf{p}\right)\right)
    \end{bmatrix} \in \mathbb{R}^{n_z}.
\end{equation}
For $\epsilon > 0$, a zero residual instead certifies a modified KKT point, whose distance from the exact KKT point grows with $\epsilon$ \cite{liao-mcphersonRegularizedSmoothedFischer2019}.
This distance of the modified KKT point from the exact KKT point does not affect the constraint satisfaction of the original problem, because even for $\epsilon > 0$, a zero residual certifies that the equality and inequality constraints are satisfied to the same tolerance as the original KKT conditions, with $\epsilon$ acting as a slack on the complementarity condition and thereby affecting the optimality gap but not the feasibility of the solution.

\subsection{Learning-Based Iterative Solver for Constrained Optimization (LISCO)}
\label{sec:bg_lisco}

\begin{figure*}[pos=!htbp]
    \centering
    \includegraphics[width=.9\textwidth]{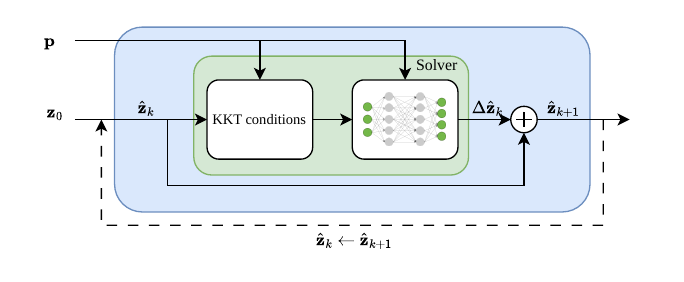}
    \caption{Monolithic LISCO update map. From the parameters $\mathbf{p}$ and the current iterate $\hat{\mathbf{z}}_k$, the modified KKT residual \eqref{eq:r_phi} is evaluated and a single dense network maps it to the full primal-dual step $\Delta\hat{\mathbf{z}}_k$, which updates the iterate. Iterating this map drives the residual to the tolerance. The decomposed map proposed in this work is shown for comparison in Figure~\ref{fig:lisco_decomp}.}
    \label{fig:lisco_monolithic}
\end{figure*}

Our previously proposed learning-based iterative solver for constrained optimization (LISCO) predicts the step $\Delta\mathbf{z}_k$ via a neural network $\Psi_{\theta}$ with parameters $\theta$~\cite{lukenSelfsupervisedLearningIterative2026,lukenLearningIterativeSolvers2024}.
At each iteration the network takes the current residual \eqref{eq:r_phi} as its input and returns the step, which is applied to the iterate $\mathbf{z}_k$ according to the update rule \eqref{eq:iterative_solver}.
We call this architecture \emph{monolithic} because a single network $\Psi_{\theta}$ maps the entire residual to the entire step.
Since every iterate $\mathbf{z}_k$ is a full primal-dual vector, the residual can be recomputed after each step and used as a stopping criterion, and each iteration requires one network evaluation with no linear-system solve.
Figure~\ref{fig:lisco_monolithic} illustrates this.

Over the course of a solve, the residual can span many orders of magnitude, which would make it a poorly scaled input.
The network therefore receives the residual split into a direction and a magnitude, the normalized residual $\mathbf{r}_k/\lVert\mathbf{r}_k\rVert_2$ and the logarithm of its norm $\log\lVert\mathbf{r}_k\rVert_2$:
\begin{equation}\label{eq:tau}
    \boldsymbol{\tau}_k = \boldsymbol{\tau}(\mathbf{z}_k,\mathbf{p})
    = \Big(\, \mathbf{r}_k / \lVert\mathbf{r}_k\rVert_2,\;\; \log \lVert\mathbf{r}_k\rVert_2 \,\Big)
    \in \mathbb{R}^{n_z + 1}.
\end{equation}
The output is rescaled by the residual norm, so that the step is large when the residual is large and shrinks to small corrections as the residual approaches zero:
\begin{equation}\label{eq:solver_step}
    \Delta\hat{\mathbf{z}}_k = \Pi_{\theta,\mathrm{solv}}(\mathbf{z}_k,\mathbf{p})
    = \gamma\, \lVert\mathbf{r}_k\rVert_2\, \Psi_{\theta}(\boldsymbol{\tau}_k,\mathbf{p}),
\end{equation}
with $\mathbf{r}_k := \mathbf{r}_\phi(\mathbf{z}_k,\mathbf{p})$ and a fixed problem-specific factor $\gamma > 0$.
The exact input scaling and the numerical safeguards are described in \cite{lukenSelfsupervisedLearningIterative2026}.

The network is trained without precomputed optimizer solutions.
For each training sample $s$, with parameter $\mathbf{p}^{(s)}$ and current iterate $\mathbf{z}_k^{(s)}$, one step produces the next iterate $\hat{\mathbf{z}}^{(s)}(\theta) = \mathbf{z}_k^{(s)} + \Pi_{\theta,\mathrm{solv}}(\mathbf{z}_k^{(s)},\mathbf{p}^{(s)})$.
The loss measures the residual at this iterate, thereby training the network to reduce the residual toward zero, where the modified KKT conditions hold.
The per-sample loss \eqref{eq:loss_per_sample} and the batch loss \eqref{eq:loss_full}, averaged over a mini-batch of $N_{\mathrm{batch}}$ parameter samples, are
\begin{equation}\label{eq:loss_per_sample}
    l(\mathbf{z},\mathbf{p}) = \log\left(\tfrac{1}{2}\lVert \mathbf{r}_{\phi}(\mathbf{z},\mathbf{p}) \rVert_2^2 + \sigma\right),
\end{equation}
\begin{equation}\label{eq:loss_full}
    L(\theta) = \frac{1}{N_{\mathrm{batch}}} \sum_{s=1}^{N_{\mathrm{batch}}} l\big(\hat{\mathbf{z}}^{(s)}(\theta),\, \mathbf{p}^{(s)}\big).
\end{equation}
The logarithm limits the variation in the loss across the orders of magnitude spanned by the residual, and the offset $\sigma > 0$ avoids numerical issues near a zero residual.
Because the loss is evaluated on the residual after the full step, the network learns to predict steps of appropriate magnitude at $\alpha_k = 1$, and the line search is kept only as a safeguard.
Local minima of the per-sample loss coincide with iterates that satisfy the modified KKT conditions, so the training target is well-posed (Lemma~1 of \cite{lukenSelfsupervisedLearningIterative2026}).
During training, the predicted steps $\Delta \mathbf{z}_k^{(s)}$ are applied to the current iterate, so each iterate $\mathbf{z}_k^{(s)}$ is generated by the solver itself.
This forward sampling trains the network on the iterates that occur at inference.
Once a sample reaches the tolerance, its parameter $\mathbf{p}^{(s)}$ is resampled and the initial iterate $\mathbf{z}_0^{(s)}$ is reset.
Section~\ref{sec:method_cl_resampling} describes the parameter resampling distribution and its closed-loop modification.

\subsection{Block-Banded Structure of the NMPC KKT System}
\label{sec:bg_banded}

The dynamics of the NMPC problem \eqref{eq:nmpc_problem} couple only consecutive stages, where each equality constraint relates $\mathbf{x}_{i+1}$ to $\mathbf{x}_i$ and $\mathbf{u}_i$, so the coupling is limited to neighboring stages along the horizon.
Additionally, if the inequality constraints are formulated as stage constraints, such as box constraints on the inputs and states, then the staged structure induces a partition of the primal-dual vector into equally sized stages:
\begin{align}\label{eq:stage_partition}
    &\mathbf{z} = \big(\mathbf{z}_0,\ldots,\mathbf{z}_{N-1}\big),
    \quad \notag\\
    &\mathbf{z}_i =
    \big(
        \underbrace{\mathbf{u}_i,\, \mathbf{x}_{i+1}}_{n_w^{s}},\;
        \underbrace{\boldsymbol{\nu}_i}_{n_h^{s}},\;
        \underbrace{\boldsymbol{\lambda}_i}_{n_g^{s}}
    \big) \in \mathbb{R}^{n_z^{s}},
\end{align}
with $n_z^{s} = n_w^{s} + n_h^{s} + n_g^{s}$ and $n_z = N\,n_z^{s}$.
The per-stage dimensions are $n_w^{s} = n_x + n_u$, $n_h^{s} = n_x$ for the state dynamics, and $n_g^{s}$ for the inequality constraints.
This partition assigns each stage $i$ exactly one input $\mathbf{u}_i$, one state $\mathbf{x}_{i+1}$, one dynamics constraint and its multiplier $\boldsymbol{\nu}_i$, together with the stage-local inequality constraints and their multipliers $\boldsymbol{\lambda}_i$.

We group the residual components according to the stage partition of the primal-dual vector:
\begin{align}\label{eq:r_phi_stages}
    &\mathbf{r}_{\phi}(\mathbf{z},\mathbf{p}) \;=\; \big(\mathbf{r}_{\phi,0},\ldots,\mathbf{r}_{\phi,N-1}\big),
    \quad \notag \\
    &\mathbf{r}_{\phi,i} =
    \begin{bmatrix}
        \nabla_{(\mathbf{u}_i,\mathbf{x}_{i+1})}\mathcal{L}\\
        \mathbf{h}_i\\
        \boldsymbol{\phi}_i
    \end{bmatrix} \in \mathbb{R}^{n_z^{s}}.
\end{align}
Each stage residual $\mathbf{r}_{\phi,i}$ has the same dimension $n_z^{s}$ as the stage block $\mathbf{z}_i$.

The discretized dynamics $\mathbf{h}_i$ involve $\mathbf{x}_i$, $\mathbf{u}_i$ and $\mathbf{x}_{i+1}$ and therefore couple stages $i$ and $i-1$ through $\mathbf{x}_i$, while the stage cost $\ell(\mathbf{x}_i,\mathbf{u}_i)$ also contains elements from stages $i$ ($\mathbf{u}_i$) and $i-1$ ($\mathbf{x}_i$).
The stationarity row of $\mathbf{x}_{i+1}$ involves $\boldsymbol{\nu}_i$ and $\boldsymbol{\nu}_{i+1}$ and so couples stage $i$ to stage $i+1$.
The stage constraints $\mathbf{g}_{\mathrm{MPC}}(\mathbf{x}_{i+1},\mathbf{u}_i)$ depend only on the current stage $i$.
Therefore, no coupling extends further than one stage, which gives the bandedness:
\begin{equation}\label{eq:banded}
    \frac{\partial\, \mathbf{r}_{\phi,i}}{\partial\, \mathbf{z}_{i'}} = \mathbf{0}
    \qquad \text{for all } \lvert i - i' \rvert > 1 .
\end{equation}

Consequently, every generalized Jacobian of $\mathbf{r}_\phi$ with respect to $\mathbf{z}$ is block-tridiagonal in the stage index, with stage bandwidth one.
Optimizers such as IPOPT \cite{wachterImplementationInteriorpointFilter2006} exploit this sparsity by utilizing sparse linear solvers such as MA27 \cite{duff1983multifrontal} to reduce the computational complexity of solving the KKT system.

While this sparse structure is inherent to the NMPC problem \eqref{eq:nmpc_problem}, the monolithic update \eqref{eq:solver_step} does not exploit it and instead couples every residual entry to every entry of the step, so the input and output dimensions of $\Psi_\theta$ grow linearly in $N$.
Under the usual choice that the hidden layers are at least as wide as the input and output dimensions, the parameter count then grows at least quadratically in $N$ \cite{lukenSelfsupervisedLearningIterative2026}.

\section{A Learned Solver for NMPC via Temporal Decomposition}\label{sec:method}

\begin{figure*}[pos=!htbp]
    \centering
    \includegraphics[width=.9\textwidth]{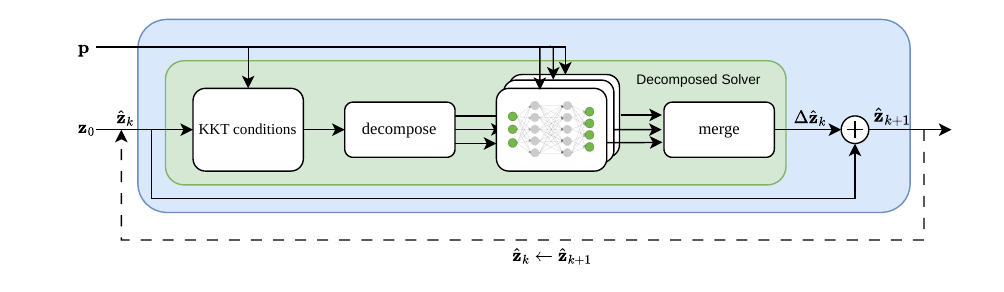}
    \caption{LISCO update map with temporal decomposition. The global modified KKT residual \eqref{eq:r_phi} is split into overlapping stage windows. A single shared network is evaluated on all windows in parallel. The overlapping predictions are combined by the fixed averaging operator. The per-stage magnitude and the global step scaling $\gamma$ are restored to form the primal-dual step.}
    \label{fig:lisco_decomp}
\end{figure*}

To exploit the block-banded structure of the NMPC problem \eqref{eq:nmpc_problem} (see Section~\ref{sec:bg_banded}), we decompose the update map along the horizon.
Figure~\ref{fig:lisco_decomp} illustrates the main idea. 
A shared neural network computes local primal-dual updates from overlapping windows of residuals from consecutive stages, and these local predictions are combined into a full update.
The network is evaluated on all windows in parallel, and predictions for stages covered by multiple windows are averaged.
Because the network is shared, its input and output sizes depend only on the window length and remain independent of the horizon.
The outer iteration \eqref{eq:iterative_solver}, the residual \eqref{eq:r_phi} and the training loss \eqref{eq:loss_full} are kept exactly as in the monolithic solver \cite{lukenSelfsupervisedLearningIterative2026}, so only the map from the residual to the step \eqref{eq:solver_step} is replaced.
The following subsections define the overlapping windows illustrated in Figure~\ref{fig:windows}, the shared network and the averaging of local predictions.

\subsection{Windows and Local Residual Features}
\label{sec:method_windows}

\begin{figure}[pos=!htbp]
    \centering
    \includegraphics[width=.9\columnwidth]{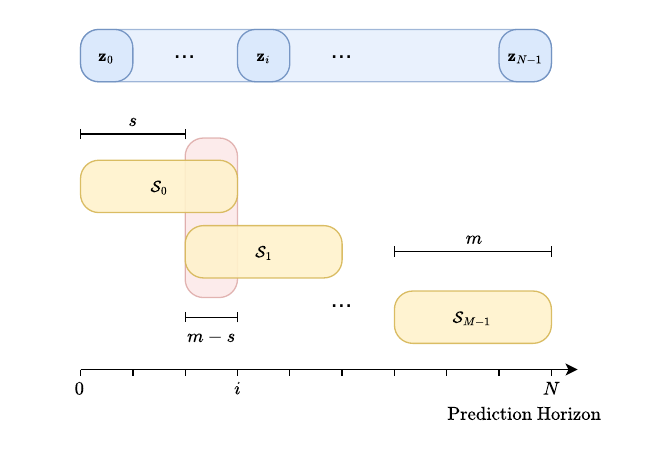}
    \caption{Time decomposition of the horizon. The stage-partitioned primal-dual iterate $\mathbf{z} = (\mathbf{z}_0,\ldots,\mathbf{z}_{N-1})$ is covered by overlapping windows $\mathcal{S}_0,\ldots,\mathcal{S}_{M-1}$ of $m$ consecutive stages, shifted by the stride $s$, so consecutive windows overlap in $m-s$ stages (red). A stage that lies in two windows, such as stage $2$, receives one prediction per window, and these are averaged into a single step by the merge operation \eqref{eq:dec_merge}. The number of windows is $M = 1 + (N-m)/s$.}
    \label{fig:windows}
\end{figure}

The proposed decomposition introduces two hyperparameters, the window length $m$ and the stride $s$.
A window is a set of $m$ consecutive stages of the horizon, and consecutive windows are shifted by $s$ stages, giving the window stage sets
\begin{align}\label{eq:dec_windows}
    &\mathcal{S}_j = \{\, js,\; js+1,\; \ldots,\; js+m-1 \,\},
    \qquad \\
    &M = 1 + \frac{N-m}{s},
\end{align}
for $j = 0,\ldots,M-1$, where $M$ is the resulting number of windows. %
We index stages by $i$, windows by $j$, and local positions within a window by $l$, while $k$ is the solver iteration.

The windows are chosen to cover the horizon with an overlap of $m-s$ stages between consecutive windows, as shown in Figure~\ref{fig:windows}.
The windows must satisfy the following admissibility conditions:
\begin{subequations}\label{eq:dec_admissibility}
    \begin{align}
        &1 \leq s \leq m \leq N,\\
        &(N-m) \bmod s = 0.
    \end{align}
\end{subequations}

The set $\mathcal{W}_i$ collects the pairs $(j,l)$ of windows $j$ and local positions $l$ inside these windows that map to stage $i$ and is defined as
\begin{equation}\label{eq:dec_owner}
    \mathcal{W}_i = \big\{\, (j,l) \;:\; js + l = i,\;\; 0 \leq j < M,\;\; 0 \leq l < m \,\big\},
\end{equation}
for all $i = 0,\ldots,N-1$.
The predictions of the windows in $\mathcal{W}_i$ are merged into a single step $\Delta \mathbf{z}_{k,i}$ for each stage $i$ by taking the average as described in Section~\ref{sec:method_update}, for the positions that overlap.

Using the stage partition \eqref{eq:r_phi_stages}, the per-stage residual $\mathbf{r}_{k,i}$, its norm $\eta_{k,i}$, and its normalization $\hat{\mathbf{r}}_{k,i}$ are given as follows:
\begin{subequations}\label{eq:dec_features}
    \begin{align}
        &\mathbf{r}_{k,i} := \mathbf{r}_{\phi,i}(\mathbf{z}_k,\mathbf{p}),
        \qquad \\
        &\eta_{k,i} := \lVert \mathbf{r}_{k,i} \rVert_2,
        \qquad \\
        &\hat{\mathbf{r}}_{k,i} := \mathbf{r}_{k,i} / \eta_{k,i},
    \end{align}
\end{subequations}
for $i = 0,\ldots,N-1$.
Unlike the input transformation of the monolithic solver \eqref{eq:tau}, which is applied once globally, the per-stage features \eqref{eq:dec_features} are computed for every stage $i$ and then collected into windows.

The stage residual $\mathbf{r}_{k,i}$ does not encode its position in the horizon.
A single shared network therefore cannot determine where a window lies in the horizon, yet the stages are not interchangeable, since the stages near the beginning and end of the horizon play different roles, particularly because of the initial condition and horizon boundary.
We therefore augment the window input with positional information that contains the stage indices of the window, so that the shared network can distinguish otherwise identical residuals occurring at different positions.
For each stage $i$ in a horizon of length $N$, the positional information is defined as:
\begin{equation}\label{eq:dec_posenc}
    \boldsymbol{\pi}_i = \Big[\; \frac{i}{N-1}, \;\; 1 - \frac{i}{N-1} \;\Big] \in \mathbb{R}^{2},
    \qquad i = 0,\ldots,N-1 ,
\end{equation}
which is normalized by $N-1$, so its range does not change with the horizon length.

The window input consists of the normalized stage residuals of the window $j$, their log-magnitudes, and the positional information of the stages at iteration $k$:
\begin{equation}\label{eq:dec_input}
    \boldsymbol{\tau}_{k,j} =
    \Big(\,
        \big(\hat{\mathbf{r}}_{k,i}\big)_{i \in \mathcal{S}_j},\;\;
        \big(\log \eta_{k,i}\big)_{i \in \mathcal{S}_j},\;\;
        \big(\boldsymbol{\pi}_{i}\big)_{i \in \mathcal{S}_j}
    \,\Big)
    \in \mathbb{R}^{m\,n_z^{s} + 3m}.
\end{equation}

The dimension of $\boldsymbol{\tau}_{k,j}$ is therefore $m\,n_z^{s} + 3m$ and is set by the window length $m$.
The problem parameter $\mathbf{p}$ is passed to every window, as in \eqref{eq:solver_step}.

\subsection{Shared Update Network and Step Averaging}
\label{sec:method_update}

The shared network $\Psi_\theta$ maps the window input $\boldsymbol{\tau}_{k,j}$ and the problem parameters $\mathbf{p}$ to a vector of predictions for the $m$ stages in window $j$:
\begin{align}\label{eq:dec_output}
    \mathbf{d}_{k,j} &= \Psi_{\theta}\big(\boldsymbol{\tau}_{k,j},\, \mathbf{p}\big) \notag \\
    &= \big(\mathbf{d}_{k,j,0},\, \ldots,\, \mathbf{d}_{k,j,m-1}\big)
    \in \mathbb{R}^{m\,n_z^{s}},
\end{align}
where the tuple denotes vector concatenation.
Each block $\mathbf{d}_{k,j,l} \in \mathbb{R}^{n_z^{s}}$ is the prediction from window $j$ for stage $i = js+l$, before averaging and rescaling.
For fixed window length $m$, stage dimension $n_z^{s}$ and parameter dimension $n_p$, the network dimensions are independent of $N$:
\begin{equation}\label{eq:dec_dims}
    n_{\mathrm{in}} = n_p + m\,n_z^{s} + 3m,
    \qquad
    n_{\mathrm{out}} = m\,n_z^{s}.
\end{equation}
If $m \ll N$, these inputs and outputs are significantly smaller than the input and output dimensions of the monolithic solver, which depend linearly on the full horizon length $N$.
The $M$ network evaluations are independent and can be executed in parallel.

For $s < m$, a stage lies in several windows, so multiple predictions are available for the same stage block and must be combined.
The overlapping predictions are merged by taking the average:
\begin{equation}\label{eq:dec_merge}
    \bar{\mathbf{d}}_{k,i} = \frac{1}{\lvert \mathcal{W}_i \rvert} \sum_{(j,l) \in \mathcal{W}_i} \mathbf{d}_{k,j,l}
    \;\in\; \mathbb{R}^{n_z^{s}},
\end{equation}
where $\lvert\mathcal{W}_i\rvert$ is the number of windows that simultaneously cover stage $i$.
The merge is applied both during training and at inference, allowing the network to learn predictions that are consistent across overlapping windows.

The resulting update for stage $i$ is therefore given by the merged prediction $\bar{\mathbf{d}}_{k,i}$, scaled by the per-stage magnitude $\eta_{k,i}$ and the global factor $\gamma$:
\begin{equation}\label{eq:dec_step}
    \Delta\hat{\mathbf{z}}_{k,i}
    = \gamma\, \min\!\big(\eta_{k,i},\, \eta_{\max}\big)\, \bar{\mathbf{d}}_{k,i},
    \qquad i = 0,\ldots,N-1.
\end{equation}
The decomposed solver scales each stage by its own magnitude $\eta_{k,i}$, in place of the single global magnitude $\lVert\mathbf{r}_k\rVert_2$ of the monolithic solver, so a stage whose residual is already small is scaled independently of one whose residual is still large.
For numerical robustness, a cap $\eta_{\max}$ bounds the per-stage magnitude.
The full update step is then given by
\begin{equation}\label{eq:dec_step_compact}
    \Delta\hat{\mathbf{z}}_{k}
    = \big(\Delta\hat{\mathbf{z}}_{k,0},\, \ldots,\, \Delta\hat{\mathbf{z}}_{k,N-1}\big).
\end{equation}

For $m = N$, the decomposed architecture differs from \eqref{eq:solver_step} only in the per-stage residual normalization and scaling, and the inclusion of the positional information.
In this case there is a single window:
\begin{equation}\label{eq:dec_mono_limit}
    m = N
    \quad\Longrightarrow\quad
    M = 1.
\end{equation}

Algorithm~\ref{alg:decomposed_solver} summarizes the decomposed solver.
Relative to the monolithic solver, lines \ref{alg:line_features} to \ref{alg:line_scale} implement the window construction, parallel network evaluation, and merging.
The shared neural network $\Psi_\theta$ is evaluated on all $M$ windows in parallel as described in line \ref{alg:line_net}, and the resulting predictions are merged into a single step $\Delta\hat{\mathbf{z}}_k$ in lines \ref{alg:line_merge} and \ref{alg:line_scale}.
The decomposition requires no communication between windows, inner iterations, or local subproblem solves.
Within an iteration the only coupling between windows is the merge \eqref{eq:dec_merge}.
Across iterations, the only coupling is the recomputation of the global residual.
The iteration runs until the residual norm falls below the tolerance $\delta$ or the iteration cap $K_{\max}$ is reached, with $\alpha_k = 1$ by default.
The line-search safeguard applied in the reported computational studies is described in Section~\ref{sec:cs_setup}.

\begin{algorithm}[t]
\caption{Solver iteration with temporal decomposition}
\label{alg:decomposed_solver}
\begin{algorithmic}[1]
\REQUIRE Parameter $\mathbf{p}$, initial iterate $\mathbf{z}_0$, trained weights $\theta$, window length $m$, stride $s$, tolerance $\delta$, cap $K_{\max}$
\STATE $k \leftarrow 0$
\WHILE{$\lVert \mathbf{r}_{\phi}(\mathbf{z}_k,\mathbf{p}) \rVert_\infty > \delta$ \textbf{and} $k < K_{\max}$}
    \STATE $\mathbf{r}_k \leftarrow \mathbf{r}_{\phi}(\mathbf{z}_k,\mathbf{p})$ \COMMENT{Eq.~\eqref{eq:r_phi}}
    \STATE Split $\mathbf{r}_k$ into stages $\mathbf{r}_{k,i}$ and form $\hat{\mathbf{r}}_{k,i}$, $\eta_{k,i}$ \label{alg:line_features} \COMMENT{Eq.~\eqref{eq:dec_features}}
    \STATE Gather window inputs $\boldsymbol{\tau}_{k,j}$ for $j = 0,\ldots,M-1$ \COMMENT{Eq.~\eqref{eq:dec_input}}
    \STATE $\mathbf{d}_{k,j} \leftarrow \Psi_{\theta}(\boldsymbol{\tau}_{k,j},\mathbf{p})$ for all $j$ in parallel \label{alg:line_net} \COMMENT{Eq.~\eqref{eq:dec_output}}
    \STATE $\bar{\mathbf{d}}_{k,i} \leftarrow$ average over $\mathcal{W}_i$ \label{alg:line_merge} \COMMENT{Eq.~\eqref{eq:dec_merge}}
    \STATE $\Delta\hat{\mathbf{z}}_{k,i} \leftarrow \gamma \min(\eta_{k,i},\eta_{\max})\, \bar{\mathbf{d}}_{k,i}$ \label{alg:line_scale} \COMMENT{Eq.~\eqref{eq:dec_step}}
    \STATE $\mathbf{z}_{k+1} \leftarrow \mathbf{z}_k + \alpha_k \Delta\hat{\mathbf{z}}_k$ \COMMENT{Eq.~\eqref{eq:iterative_solver}}
    \STATE $k \leftarrow k+1$
\ENDWHILE
\RETURN $\mathbf{z}_k$
\end{algorithmic}
\end{algorithm}

\subsection{Self-supervised Training}
\label{sec:method_training}

We train the shared network using the self-supervised loss from our previous work \cite{lukenSelfsupervisedLearningIterative2026}, summarized in Section~\ref{sec:bg_lisco}.
The loss \eqref{eq:loss_full} measures the residual of the full problem after the decomposed update and is averaged over a batch of primal-dual iterates.
The residual and stopping criterion retain their interpretation from the monolithic solver, since a zero residual still certifies a modified KKT point.
Training iterates are generated by the forward sampling strategy of Section~\ref{sec:bg_lisco}, which uses the solver's own iterates to cover a broad range of residual magnitudes.

\subsection{Closed-Loop Resampling}
\label{sec:method_cl_resampling}

The forward sampling strategy of Section~\ref{sec:bg_lisco} generates training iterates by applying the learned solver to sampled problems.
In the original training procedure, the parameters defining these problems are resampled uniformly.
For the NMPC formulation considered here, the problem parameter is the initial state, $\mathbf{p} = \mathbf{x}_{\mathrm{init}}$.
During closed-loop operation, these initial states are the states visited by the controller.
In economic NMPC in particular, the controller may drive the system to a steady state on the boundary of the feasible set, and such states with several simultaneously active constraints are rarely produced by uniform sampling over the state space.
A solver trained on uniformly sampled initial states therefore encounters a distribution shift in closed-loop operation.

To align the training distribution, we modify the training resampling strategy by introducing a reset probability $q_{\mathrm{reset}} \in (0,1]$, as summarized in Algorithm~\ref{alg:solver_training}, which extends the self-supervised training procedure of our previous work \cite{lukenSelfsupervisedLearningIterative2026}.
Whenever a sample reaches the training tolerance $\delta_{\mathrm{train}}$, the initial state of the next training problem is set to the predicted next state $\mathbf{x}_1$ of the converged solution with probability $1 - q_{\mathrm{reset}}$.
With probability $q_{\mathrm{reset}}$, a new initial state is sampled uniformly instead.
Since $\mathbf{x}_1$ is a component of the primal solution $\hat{\mathbf{w}}$, this state is available in the converged iterate at no additional cost, and the expected trajectory length is $1/q_{\mathrm{reset}}$ steps.
For $q_{\mathrm{reset}} = 1$, every converged sample is replaced by a uniformly resampled parameter, which recovers the uniform training distribution from our previous work \cite{lukenSelfsupervisedLearningIterative2026}.

Two safeguards bound the effort spent on a single sample, the iteration cap $N_{\max}$ on the number of steps and the stall cap $N_{\mathrm{stall}}$ on the number of steps without improvement of the residual.
A sample that reaches either cap is treated as non-converged and resampled uniformly.

\begin{algorithm}[t]
\caption{Self-Supervised Solver Training with Closed-Loop Resampling ($q_{\mathrm{reset}} = 1$ recovers uniform resampling)}
\label{alg:solver_training}
\begin{algorithmic}[1]
\REQUIRE Batch size $N_{\mathrm{batch}}$, tolerance $\delta_{\mathrm{train}}$, iteration cap $N_{\max}$, stall cap $N_{\mathrm{stall}}$, reset probability $q_{\mathrm{reset}} \in (0,1]$
\STATE Sample $\{\mathbf{p}^{(s)}\}_{s=1}^{N_{\mathrm{batch}}}$, initialize the iterates $\mathbf{z}^{(s)}$, set $k^{(s)}, c^{(s)} \leftarrow 0$
\FOR{each training step}
    \STATE $\hat{\mathbf{z}}^{(s)} \leftarrow \mathbf{z}^{(s)} + \Pi_{\theta,\mathrm{solv}}(\mathbf{z}^{(s)},\mathbf{p}^{(s)})$ for all $s$ \COMMENT{Eq.~\eqref{eq:iterative_solver}, $\alpha_k = 1$}
    \STATE Update $\theta$ by a stochastic gradient step on $L(\theta)$ \COMMENT{Eq.~\eqref{eq:loss_full}}
    \FOR{each sample $s$}
        \STATE $\mathbf{z}^{(s)} \leftarrow \hat{\mathbf{z}}^{(s)}$, \quad $k^{(s)} \leftarrow k^{(s)} + 1$
        \STATE $c^{(s)} \leftarrow 0$ if $\lVert \mathbf{r}_{\phi}(\mathbf{z}^{(s)},\mathbf{p}^{(s)}) \rVert_2^2$ improves on its best value so far, else $c^{(s)} \leftarrow c^{(s)} + 1$
        \IF{$\tfrac{1}{2}\lVert \mathbf{r}_{\phi}(\mathbf{z}^{(s)},\mathbf{p}^{(s)}) \rVert_2^2 < \delta_{\mathrm{train}}$}
            \STATE $\mathbf{p}^{(s)} \leftarrow \mathbf{x}_1^{(s)}$ with probability $1 - q_{\mathrm{reset}}$, otherwise resample $\mathbf{p}^{(s)}$ \COMMENT{$\mathbf{x}_1^{(s)}$ from $\hat{\mathbf{w}}^{(s)}$}
            \STATE Reset $\mathbf{z}^{(s)}$ to the initial iterate, \quad $k^{(s)}, c^{(s)} \leftarrow 0$
        \ELSIF{$k^{(s)} \geq N_{\max}$ \textbf{or} $c^{(s)} \geq N_{\mathrm{stall}}$}
            \STATE Resample $\mathbf{p}^{(s)}$, reset $\mathbf{z}^{(s)}$ to the initial iterate, \quad $k^{(s)}, c^{(s)} \leftarrow 0$
        \ENDIF
    \ENDFOR
\ENDFOR
\RETURN trained solver $\Pi_{\theta,\mathrm{solv}}$
\end{algorithmic}
\end{algorithm}

\section{Case Study}\label{sec:case_study}

We evaluate the proposed method on the economic NMPC of a cascade of continuous stirred-tank reactors (CSTR).
The studies investigate horizon and system-size scaling, smoothing, window length and overlap, and closed-loop control performance.
We compare against monolithic LISCO, IPOPT \cite{wachterImplementationInteriorpointFilter2006} with MA27 \cite{duff1983multifrontal}, and, in closed loop, approximate MPC.

\subsection{CSTR Cascade Process Model}
\label{sec:cs_model}

We consider a cascade of $n_r$ continuous stirred-tank reactors as presented in \cite{heinleinRobustModelPredictive2025}.
A solvent stream enters the first reactor and passes through the cascade, and the reactants A and B are fed into every reactor separately.
Two reactions take place.
The first reaction, $\mathrm{A} + \mathrm{B} \rightarrow \mathrm{R}$, produces the desired product R, whereas the second reaction, $2\,\mathrm{A} \rightarrow \mathrm{S}$, produces the undesired side product S.
The objective is to maximize the production of R, minimize the production of S, limit control effort, and satisfy constraints on the states and inputs.

The states of each reactor $n$ are the concentrations of the four species, $c_{\mathrm{A},n}$, $c_{\mathrm{B},n}$, $c_{\mathrm{R},n}$, $c_{\mathrm{S},n}$ in $\mathrm{mol}/\mathrm{L}$, and the reactor temperature $T_{r,n}$ in $^\circ\mathrm{C}$.
The inputs are the molar feed rates $u_{\mathrm{A},n}$ and $u_{\mathrm{B},n}$ of A and B, respectively, in $\mathrm{mol}/\mathrm{s}$.
All reactors have equal volume $V_n = V_R/n_r$, so the total reactor volume is independent of the cascade length.
We consider a fixed jacket temperature $T_j$ shared by all reactors.
Stacking the per-reactor states and inputs gives the following vectors and dimensions:
\begin{align}\label{eq:cs_dims}
    \mathbf{x} &= \big( c_{\mathrm{A},1:n_r},\; c_{\mathrm{B},1:n_r},\; c_{\mathrm{R},1:n_r},\; c_{\mathrm{S},1:n_r},\; T_{r,1:n_r} \big) \in \mathbb{R}^{n_x}, \notag\\
    \mathbf{u} &= \big( u_{\mathrm{A},1:n_r},\; u_{\mathrm{B},1:n_r} \big) \in \mathbb{R}^{n_u}, \notag\\
    n_x &= 5\,n_r, \qquad n_u = 2\,n_r, \qquad V_n = \frac{V_R}{n_r}.
\end{align}

The reaction rates $r_{1,n}$ and $r_{2,n}$ follow Arrhenius kinetics, with pre-exponential constants $k_1$, $k_2$ and activation energies $E_{a,1}$, $E_{a,2}$:
\begin{align}\label{eq:cs_rates}
    r_{1,n} &= k_1 \exp\!\left( \frac{-E_{a,1}}{R_g\,(T_{r,n}+273.15)} \right) c_{\mathrm{A},n}\, c_{\mathrm{B},n}, \notag\\
    r_{2,n} &= k_2 \exp\!\left( \frac{-E_{a,2}}{R_g\,(T_{r,n}+273.15)} \right) c_{\mathrm{A},n}^{2}.
\end{align}
The product reaction is second order overall, first order in each reactant, while the side reaction is second order in A alone.

The ordinary differential equations (ODEs) for the $n$-th reactor are given as follows:
\begin{subequations}\label{eq:cs_ode}
    \begin{align}
        \dot c_{\mathrm{A},n} &= \beta\,(c_{\mathrm{A},n-1}-c_{\mathrm{A},n}) - r_{1,n} - 2\,r_{2,n} + \frac{u_{\mathrm{A},n}}{V_n},\\[2pt]
        \dot c_{\mathrm{B},n} &= \beta\,(c_{\mathrm{B},n-1}-c_{\mathrm{B},n}) - r_{1,n} + \frac{u_{\mathrm{B},n}}{V_n},\\[2pt]
        \dot c_{\mathrm{R},n} &= \beta\,(c_{\mathrm{R},n-1}-c_{\mathrm{R},n}) + r_{1,n},\\[2pt]
        \dot c_{\mathrm{S},n} &= \beta\,(c_{\mathrm{S},n-1}-c_{\mathrm{S},n}) + r_{2,n},\\[2pt]
        \dot T_{r,n} &= \beta\,(T_{r,n-1}-T_{r,n})
            + \frac{k_A}{\rho\,c_p\,V_n}\,(T_j - T_{r,n}) \notag\\
            &\quad + \frac{1}{\rho\,c_p}\big( -\Delta H_1\, r_{1,n} - \Delta H_2\, r_{2,n} \big),
    \end{align}
\end{subequations}
for $n = 1,\ldots,n_r$.

We impose the cascade boundary conditions $c_{\mathrm{A},0} = c_{\mathrm{B},0} = c_{\mathrm{R},0} = c_{\mathrm{S},0} = 0$ and $T_{r,0} = T_{r,\mathrm{in}}$, with the dilution rate $\beta = V_{\mathrm{out}}/V_n$.
Thus, the stream entering the first reactor contains no reactants or products and has the inlet temperature, so the same balance equations apply to the first and all downstream reactors.
The balances in \eqref{eq:cs_ode} account for convection from the upstream reactor, outflow from the current reactor, reaction production and consumption, and the feeds of A and B.
The temperature balance additionally includes jacket heat transfer and reaction enthalpies.
Because the balance for reactor $n$ depends only on reactors $n-1$ and $n$, the cascade has local spatial coupling.
This spatial structure is not exploited in the present work and is left for future investigation.
Stacking the balances of all reactors gives the continuous-time model $\dot{\mathbf{x}} = f_c(\mathbf{x},\mathbf{u})$.

The model is discretized using an implicit Euler scheme with step $\Delta t$ as follows:
\begin{equation}\label{eq:cs_euler}
    \mathbf{x}_{i+1} = \mathbf{x}_i + \Delta t\; f_c(\mathbf{x}_{i+1},\mathbf{u}_i) =: f(\mathbf{x}_i,\mathbf{u}_i),
\end{equation}
which defines the discrete dynamics $f$ used in \eqref{eq:nmpc_problem} with time and stage index $i$.
We use one implicit-Euler step per NMPC interval.
The resulting discrete state equations form equality constraints of the NLP and preserve the stage structure of Section~\ref{sec:bg_nmpc}.
In this case study, we use $\Delta t = 1\,\mathrm{s}$.
The physical parameters are given in Table~\ref{tab:cstr_params}. %
The reactors are cooled, since the jacket temperature is below the inlet temperature, and both reactions are exothermic.

\begin{table}[pos=!htbp]
\centering
\caption{Parameters of the CSTR cascade model. The first reaction $\mathrm{A}+\mathrm{B} \rightarrow \mathrm{R}$ is denoted by subscript $1$ and the second reaction $2\,\mathrm{A} \rightarrow \mathrm{S}$ by subscript $2$.}
\label{tab:cstr_params}
\begin{tabular}{lll}
\toprule
Symbol & Value [Unit] & Description \\
\midrule
$n_r$              & study dependent [--]        & Number of reactors \\
$\Delta t$         & 1.0    [s]                  & Discretization step \\
$V_R$              & 15     [L]                  & Total reactor volume \\
$V_{\mathrm{out}}$ & 1      [L/s]                & Volumetric flow rate \\
$k_1$              & 2      [L/(mol\,s)]         & Rate constant \\
$k_2$              & 2      [L/(mol\,s)]         & Rate constant \\
$E_{a,1}$          & 500    [J/mol]              & Activation energy \\
$E_{a,2}$          & 600    [J/mol]              & Activation energy \\
$R_g$              & 8.3145 [J/(mol\,K)]         & Universal gas constant \\
$\Delta H_1$       & $-100$ [kJ/mol]             & Reaction enthalpy \\
$\Delta H_2$       & $-50$  [kJ/mol]             & Reaction enthalpy \\
$k_A$              & 2      [kW/K]               & Heat transfer coeff. $\times$ area \\
$\rho$             & 1      [kg/L]               & Density \\
$c_p$              & 0.5    [kJ/(kg\,K)]         & Heat capacity \\
$T_j$              & 50     [$^\circ$C]         & Jacket temperature \\
$T_{r,\mathrm{in}}$& 60     [$^\circ$C]        & Inlet temperature \\
\bottomrule
\end{tabular}
\end{table}

\subsection{Economic NMPC Formulation}
\label{sec:cs_nmpc}

The economic NMPC considered here is an instance of the NMPC problem \eqref{eq:nmpc_problem} with the discrete model \eqref{eq:cs_euler}, without a terminal cost and terminal constraints.
The parameter is the current state, so $n_p = n_x$, and the initial state is eliminated from the decision vector and enters only the first dynamics constraint through $\mathbf{x}_0 = \mathbf{p}$.

We use the following economic stage cost, with the weights $w_R=10$, $w_S=1$, $w_{u_A}=0.1$ and $w_{u_B}=0$:
\begin{align}\label{eq:cs_objective}
    \ell(\mathbf{x}_{i+1},\mathbf{u}_i) =
    &- w_R\, c_{\mathrm{R},n_r}(\mathbf{x}_{i+1}) \notag \\
    &+ w_S\, c_{\mathrm{S},n_r}(\mathbf{x}_{i+1}) \notag \\
    &+ w_{u_A} \lVert \mathbf{u}_{\mathrm{A},i} \rVert_2^2 \notag \\
    &+ w_{u_B} \lVert \mathbf{u}_{\mathrm{B},i} \rVert_2^2 ,
\end{align}
where $c_{\mathrm{R},n_r}$ and $c_{\mathrm{S},n_r}$ are the product and side-product concentrations of the last reactor and $\mathbf{u}_{\mathrm{A},i}$, $\mathbf{u}_{\mathrm{B},i}$ are the feed rates of all reactors at stage $i$.
Only the concentrations in the last reactor enter the product-related terms because this reactor defines the cascade outlet.
The feed of B is not penalized, $w_{u_B} = 0$, and is limited only by its bounds.
Since the objective is economic with strong weights on the maximization of the product and minimization of the side product, the NMPC drives the system to a steady state on the boundary of the feasible set.

Each dynamics constraint $i$ couples $\mathbf{x}_i$, $\mathbf{u}_i$ and $\mathbf{x}_{i+1}$, which is the nearest-neighbor coupling of the stage partition \eqref{eq:stage_partition}.
The inequality constraints $\mathbf{g}_{\mathrm{MPC}}(\mathbf{x}_{i+1},\mathbf{u}_i) \leq \mathbf{0}$ in \eqref{eq:nmpc_problem} are formulated as box constraints on the states and inputs.
Since the physical variables differ by orders of magnitude, from concentrations of order $10^{-1}\,\mathrm{mol}/\mathrm{L}$ to temperatures of order $10^{2}\,^\circ\mathrm{C}$, the resulting NLP may be ill-conditioned.
To improve the conditioning, the decision variables are scaled so that their typical values lie roughly between 0 and 1.
The scaling ranges and the bounds of the decision variables are given in Table~\ref{tab:cstr_scaling}.

\begin{table}[pos=!htbp]
\centering
\caption{Scaling ranges and decision variable bounds. The listed scaling ranges are mapped to $[0,1]$.}
\label{tab:cstr_scaling}
\begin{tabular}{llll}
\toprule
Variable & Unit & Scaling range & Bounds \\
\midrule
$c_{\mathrm{A},n}$ & mol/L      & $[0,\;0.1]$  & $[0,\;\infty)$ \\
$c_{\mathrm{B},n}$ & mol/L      & $[0,\;1]$    & unconstrained \\
$c_{\mathrm{R},n}$ & mol/L      & $[0,\;4]$    & unconstrained \\
$c_{\mathrm{S},n}$ & mol/L      & $[0,\;1]$    & unconstrained \\
$T_{r,n}$          & $^\circ$C  & $[20,\;80]$  & $[20,\;80]$ \\
$u_{\mathrm{A},n}$ & mol/s      & $[0,\;1.5]$  & $[0,\;1.5]$ \\
$u_{\mathrm{B},n}$ & mol/s      & $[0,\;1.5]$  & $[0,\;1.5]$ \\
\bottomrule
\end{tabular}
\end{table}

The resulting per-stage dimensions of the stage partition \eqref{eq:stage_partition} as a function of the number of reactors $n_r$ are
\begin{align}\label{eq:cs_nlp_dims}
    n_w^{s} &= 7\,n_r   && \text{primal variables,} \notag \\
    n_h^{s} &= 5\,n_r   && \text{equality constraints,} \notag \\
    n_g^{s} &= 7\,n_r   && \text{inequality constraints,} \notag \\
    n_z^{s} &= 19\,n_r  && \text{primal-dual variables,} \notag \\
    n_p     &= 5\,n_r   && \text{parameters.}
\end{align}

\subsection{Computational Studies}
\label{sec:cs_baselines}

The decomposed solver is evaluated in four computational studies, which investigate the scaling with the horizon and the number of reactors, the effect of the smoothing parameter $\epsilon$, the influence of the window length and stride, and the closed-loop control performance.
Throughout these studies, four methods for solving the NMPC problem \eqref{eq:nmpc_problem} are compared.
Monolithic LISCO is the architecture presented in our previous work \cite{lukenSelfsupervisedLearningIterative2026}.
The decomposed solver of Section~\ref{sec:method} is the proposed method and we use IPOPT \cite{wachterImplementationInteriorpointFilter2006} as the NLP-solver baseline.
Approximate MPC maps the initial state directly to the first control action, is trained by imitation learning on precomputed optimizer solutions, and is considered only in the closed-loop comparison.
The code for reproducing all studies is openly available.\footnote{\url{https://github.com/lukaslueken/lisco-decompose}}

All controllers are based on the same NMPC formulation and are evaluated on the same parameter samples, so speedups are compared per problem instance.
The two learned solvers use the same hidden architecture, consisting of two hidden layers of 1024 units with GELU activation \cite{hendrycksGaussianErrorLinear2016}, and share the same training budget, optimizer settings, loss \eqref{eq:loss_full}, resampling strategy, and evaluation settings.
Their input and output layers necessarily differ, since by \eqref{eq:solver_step} the monolithic network maps the full residual to the full step, so these layers scale with $N\,n_z^{s}$, whereas by \eqref{eq:dec_dims} the decomposed network processes a single window and is independent of $N$.
The monolithic network therefore carries more trainable parameters than the decomposed one, with the gap widening with $N$ and $n_r$.
Parameter counts and training times are reported per configuration in Section~\ref{sec:res_scaling}.

IPOPT is configured with the sparse symmetric-indefinite linear solver MA27 from the HSL library \cite{duff1983multifrontal}, which exploits the block-banded structure of \eqref{eq:banded}, at a tolerance of $10^{-5}$.
Comparing solve times between the learned solvers and IPOPT further requires that both target a comparable accuracy.
The learned solver runs at a fixed smoothing $\epsilon$, which relaxes the complementarity condition by up to $\epsilon^2/2$, and the interior-point barrier parameter plays the same role for IPOPT, so the barrier target is set to $\mu_{\mathrm{target}} = \epsilon^2/2$.
With $\epsilon = 10^{-1}$ this gives $\mu_{\mathrm{target}} = 5\cdot 10^{-3}$, which is used in the horizon-scaling and closed-loop studies.

In all studies, the solvers are cold-started to directly measure the quality of the predicted steps rather than the quality of an initial guess.
The learned solvers initialize the primal variables at the midpoint of the scaling ranges and the multipliers at zero, in training and at inference, and IPOPT starts from its default initialization.

The horizon-scaling study evaluates how solver performance scales with the horizon length $N$ and the number of reactors $n_r$.
Monolithic LISCO and the decomposed solver are trained and then evaluated alongside IPOPT for the horizon lengths $N \in \{10, 20, 50, 100\}$ and the reactor counts $n_r \in \{2, 3, 5\}$, with all remaining settings fixed by Table~\ref{tab:setup}.
Training effort, parameter counts, convergence and solve times are reported in Section~\ref{sec:res_scaling}.

The choice of the smoothing parameter $\epsilon$ involves a trade-off, as larger values improve the conditioning of the residual but relax the complementarity condition further.
This trade-off is quantified in the smoothing study, in which the decomposed solver is trained and evaluated for a sweep of $\epsilon$ at a fixed configuration with $n_r = 2$ and $N = 50$.
The resulting optimality gap is evaluated in Section~\ref{sec:res_eps} against IPOPT, which serves here as the objective reference and is therefore run at a tolerance of $10^{-8}$ without an increased barrier parameter.

In the overlap study, we vary the window length $m$ and the stride $s$ of the decomposition, which determine the locality of the update map and the number of predictions per stage.
This study considers only the decomposed solver, with all remaining settings taken from Table~\ref{tab:setup}, and the results are reported in Section~\ref{sec:res_overlap}.

In the closed-loop study, we apply the decomposed solver as an NMPC controller and compare it to IPOPT and approximate MPC.
The three controllers start from the same set of 50 uniformly sampled initial states at a fixed configuration with $n_r = 5$ and $N = 50$.
Each trajectory is simulated for 50 steps, and the closed-loop cost and constraint violations are compared in Section~\ref{sec:res_closed_loop}.
In this study, the decomposed solver is trained with the closed-loop resampling of Section~\ref{sec:method_cl_resampling} at $q_{\mathrm{reset}} = 0.02$, so the expected trajectory length of 50 steps during training matches the simulation length.
All remaining studies use uniform resampling with $q_{\mathrm{reset}} = 1$.
To isolate the effect of the training distribution, a second solver trained with uniform resampling and otherwise identical settings is evaluated in the same closed-loop setting.
Approximate MPC, as introduced in Section~\ref{sec:bg_nmpc}, is trained on closed-loop data.
We collect 150000 optimizer solutions from 3000 closed-loop trajectories with uniformly sampled initial states and 50 steps each.
These solutions are generated with IPOPT using a tolerance of $10^{-5}$ and a barrier target of $\mu_{\mathrm{target}} = 5\cdot10^{-3}$ to match the accuracy of the learned solvers and the IPOPT baseline.
Its training uses the same hidden network architecture as the learned solvers, a batch size of 1024, a learning rate of $5\cdot10^{-4}$ and 2500 epochs.
Generating its training set requires solving the NLP for every sample, sequentially along each trajectory, whereas the loss \eqref{eq:loss_full} requires no reference solutions.
Section~\ref{sec:res_closed_loop} compares the resulting offline computational effort.

\subsection{Training, Evaluation and Implementation}
\label{sec:cs_setup}

All four studies share one configuration, which is listed in Table~\ref{tab:setup}.
The overlap and closed-loop studies use $n_r = 5$ reactors, while the smoothing study uses $n_r = 2$.
The horizon-scaling study varies the number of reactors.
Following Algorithm~\ref{alg:solver_training}, the training loop evaluates the loss after applying one solver step, then applies these steps to the current iterates and resamples the parameters of samples that reach the tolerance $\delta_{\mathrm{train}}$, the iteration cap $N_{\max}$ or the stall cap $N_{\mathrm{stall}}$.
The number of epochs scales with the horizon as $N_{\mathrm{epochs}} = 6000\,N$, so that for longer horizons, training takes proportionally longer, with the same rule for both learned architectures.
The default window length and stride, $m = 5$ and $s = 1$, are used in all studies except the overlap study, which examines the influence of this choice.
The NLP and its residuals are evaluated in double precision, while the network operates in bfloat16.
This halves the network memory relative to single precision without a meaningful loss in accuracy, since the network output is scaled relative to the residual magnitude and the residual itself is evaluated in double precision.
All runs use $\Delta t = 1\,\mathrm{s}$, so the horizon length and the predicted time interval scale together. %

\begin{table}[pos=!htbp]
\centering
\caption{Shared configuration of all studies.}
\label{tab:setup}
\begin{tabular}{lll}
\toprule
Group & Setting & Value \\
\midrule
Network      & hidden layers          & $2 \times 1024$, GELU \cite{hendrycksGaussianErrorLinear2016} \\
             & precision              & bfloat16 \\
\midrule
Solver       & step scaling $\gamma$  & 0.1 \\
             & norm clamp $\eta_{\max}$ & 100 \\
             & smoothing $\epsilon$   & $10^{-1}$ \\
             & penalization $\rho_{\mathrm{FB}}$ & 0.8 \\
             & window $m$ / stride $s$ & 5 / 1 \\
\midrule
Training     & epochs                 & $6000\,N$ \\
             & batch size             & 256 \\
             & optimizer              & AdamW \cite{loshchilov2017decoupled} \\
             & learning rate          & $5\cdot10^{-4}$ \\
             & regularization         & $\rho_{\mathrm{reg}} = \omega = 10^{-2}$ \\
             & loss offset $\sigma$ & $10^{-16}$ \\
             & reset prob. $q_{\mathrm{reset}}$ & 1 (0.02 closed loop) \\
             & tolerance $\delta_{\mathrm{train}}$ & $\tfrac{1}{2}\lVert \mathbf{r}_\phi \rVert_2^2 < 10^{-16}$ \\
             & stall cap $N_{\mathrm{stall}}$ & 2 \\
             & iteration cap $N_{\max}$ & 1000 \\
\midrule
Evaluation   & samples                & study specific \\
             & tolerance $\delta$     & $\lVert \mathbf{r}_\phi \rVert_\infty \leq 10^{-5}$ \\
             & iteration cap          & study specific \\
             & initialization         & cold start \\
\bottomrule
\end{tabular}
\end{table}

While the objective function \eqref{eq:cs_objective} is convex, the nonlinear dynamics \eqref{eq:cs_euler} cause the NLP to be nonconvex.
Therefore, the Hessian of the Lagrangian can be indefinite, and Assumption~5 of \cite{lukenSelfsupervisedLearningIterative2026}, which underlies the well-posedness result of Section~\ref{sec:bg_lisco}, may not hold.
To keep the training problem well-posed, we form the training gradient with the convexification and regularization of \cite{lukenSelfsupervisedLearningIterative2026}.
We replace the Hessian of the Lagrangian with the objective Hessian, which is diagonal with entries $2\,w_{u_A}$ on the feed rates of A and zero elsewhere, and add to its diagonal $\rho_{\mathrm{reg}}\mathbf{I}$ with the regularization factor $\rho_{\mathrm{reg}} > 0$ so that it is positive definite.
We further apply the same dual regularization as in \cite{lukenSelfsupervisedLearningIterative2026} with factor $\omega > 0$ to improve the conditioning of the KKT matrix during training when the multipliers grow large.
Both modifications act only on the gradient used to update the network during training and leave the residual \eqref{eq:r_phi} unchanged, so convergence at inference and all reported results are those of the original problem.

For the evaluation, a problem instance is drawn by uniformly sampling the initial state from the scaling ranges listed in Table~\ref{tab:cstr_scaling} and retaining it only if IPOPT finds a feasible solution, so that all controllers are compared on instances known to be feasible.
The iteration cap is 300 in the horizon-scaling, smoothing and closed-loop studies and 200 in the overlap study.
A sample that does not reach the tolerance within the iteration cap is counted as non-convergent rather than discarded, so the convergence fraction and the timing statistics refer to the same sample set.
For a non-convergent sample, the solver returns the iterate with the lowest residual norm encountered, whose first input is applied in the closed-loop simulation.
The reported solve times include all necessary computations for its implementation: residual evaluation, feature construction, network evaluation, merging, scaling and the iterate update, over all iterations until the tolerance is met.
The one-time setup and compilation times are excluded.
Timings are measured on single instances without batching, because an NMPC controller solves one problem per sampling instant.
The learned solvers apply their steps with a two-way backtracking line search.
The default step size is $\alpha_k = 1$, and it is reduced by a factor of $0.95$ only if the residual norm of an iterate exceeds the best value reached so far by at least a factor of $10$.
Once the residual decreases again, the step size is increased by a factor of $1/0.95$ per iteration, up to a maximum of $1$.
Since the solver is trained to predict steps of appropriate size, the line search is activated only for rare, ill-conditioned instances.

All computational studies are run on an AMD Ryzen 9 9950X with 16 cores and an NVIDIA GeForce RTX 4090 under WSL2.
The NLP, its derivatives and the modified KKT residual are implemented in CasADi \cite{anderssonCasADiSoftwareFramework2019}, the networks in PyTorch \cite{paszkePyTorchImperativeStyle2019}, and during training the CasADi functions are evaluated in parallel over the batch on the GPU via cusadi \cite{jeonCusADiGPUParallelization2024}, so the GPU is used only during training.
At inference, all computations run on the CPU, where the KKT functions are executed as C code generated by CasADi, the solver iteration together with the decomposition and merge logic is compiled with Numba \cite{lamNumbaLLVMbasedPython2015}, and the network is evaluated in bfloat16 with PyTorch.
The learned solvers and IPOPT are thus compared on the same hardware.

\section{Results}\label{sec:results}

\subsection{Horizon and Problem-Size Scaling}
\label{sec:res_scaling}

\begin{figure}[pos=!htbp]
\centering
\includegraphics[width=.9\columnwidth]{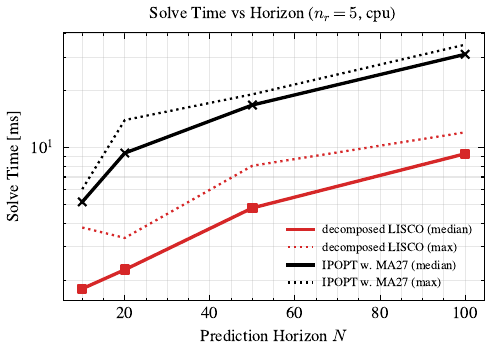}
\caption{CPU solve time against the horizon length $N$ at $n_r = 5$ for the decomposed solver and IPOPT with MA27, with the median (solid) and the maximum (dotted) over 256 instances per configuration. The monolithic solver is not shown, as it does not converge at $n_r = 5$.}
\label{fig:decomposed_solver_vs_monolithic_solver_vs_ipopt}
\end{figure}

For the horizon-scaling study, we train the monolithic and decomposed solvers for the horizon lengths $N \in \{10, 20, 50, 100\}$ and the reactor counts $n_r \in \{2, 3, 5\}$ and evaluate them against IPOPT, with all remaining settings fixed by Table~\ref{tab:setup}.
Table~\ref{tab:horizon_training} lists the convergence fractions, the trainable parameter counts and the training times of the learned solvers, and the solve times are shown in Fig.~\ref{fig:decomposed_solver_vs_monolithic_solver_vs_ipopt} and Table~\ref{tab:horizon_runtime}.

The proposed decomposed solver reaches the residual tolerance on all 256 evaluation samples for every combination of horizon length and reactor count.
The monolithic solver converges on all samples only at $n_r = 2$ and $N = 10$, on 94.5\% and 99.6\% of the samples at $(n_r, N) = (2, 20)$ and $(3, 10)$, and on no sample in the remaining configurations.
No monolithic runs are conducted at $N = 100$, since the architecture already fails to converge at shorter horizons under the given network size.
Both learned solvers share the hidden architecture and the training budget of Table~\ref{tab:setup}, and the monolithic network has up to 5.2 times as many trainable parameters, with 10.81 versus 2.06 million parameters at $n_r = 5$ and $N = 50$, the largest configuration with a monolithic run, because its input and output layers grow with $N\,n_z^{s}$ by \eqref{eq:solver_step}.
To restore convergence, the monolithic network and its training effort would have to grow with the problem size \cite{lukenSelfsupervisedLearningIterative2026}.
The decomposition avoids this cost, as the input and output dimensions of the shared network are set by the window length through \eqref{eq:dec_dims}.
Since the number of epochs scales with $N$, the training times grow with the horizon and reach 2.81 hours for the decomposed solver at $n_r = 5$ and $N = 100$.
For the same number of training epochs, training the decomposed solver takes 1.16 to 1.75 times longer than training the monolithic solver, since the shared network is evaluated $M$ times per step, once per window, which outweighs its smaller size.

Fig.~\ref{fig:decomposed_solver_vs_monolithic_solver_vs_ipopt} shows the CPU solve times over the horizon length at $n_r = 5$, and Table~\ref{tab:horizon_runtime} lists the timing statistics for all configurations.
All solvers run on the same CPU, and IPOPT is configured for a comparable relaxation using the barrier target of Section~\ref{sec:cs_baselines}.
Under these conditions the decomposed solver achieves lower median solve times than IPOPT with MA27 in every configuration, with median speedups between 1.68 and 4.11, the largest occurring at $n_r = 5$ and $N = 20$.
The advantage also holds for the upper tail of the solve-time distribution, as the 95th percentile solve time of the decomposed solver lies below the median solve time of IPOPT in every configuration, with similar tail behavior across the reactor counts.
The per-sample speedup remains above one in every configuration except at $(n_r, N) = (2, 20)$ and $(2, 100)$, where its minimum is 0.52 and 1.00, respectively.
In each iteration of the decomposed solver, the shared network is evaluated in parallel for $M$ windows, and the reported times include this cost.
The measured times remain specific to the hardware and implementation of Section~\ref{sec:cs_setup}.

We use the same window length, stride and hidden architecture at every reactor count, without retuning the decomposition.
The size of the decomposed network is independent of the horizon length and grows only mildly with the reactor count through the stage dimension $n_z^{s} = 19\,n_r$ of \eqref{eq:cs_nlp_dims}.
The convergence of the decomposed solver is unaffected by the reactor count, and its median solve times increase from $3.44\,\mathrm{ms}$ to $4.81\,\mathrm{ms}$ at $N = 50$.
As shown in Table~\ref{tab:horizon_runtime}, the median speedup over IPOPT increases with the reactor count for every tested horizon length $N \geq 20$.
Further speedups are possible if a good initial guess is provided, e.g. via a predictor strategy as in \cite{lukenSelfsupervisedLearningIterative2026}.

\begin{table*}[pos=!htbp]
\centering
\caption{Training cost and convergence of the monolithic and the decomposed solver at equal hidden architecture and training budget for all reactor counts $n_r$ and horizon lengths $N$. Reported are the trainable parameters $n_\mathrm{par}$, the training times $t_\mathrm{train}$ and the fraction of the 256 evaluation samples reaching $\lVert \mathbf{r}_\phi \rVert_\infty \leq 10^{-5}$ within 300 iterations, where 0.000 means that no sample reached the tolerance. Dashes mark configurations without a monolithic run.}
\label{tab:horizon_training}
\begin{tabular}{ccrrrrrr}
\toprule
 &  & \multicolumn{3}{c}{monolithic} & \multicolumn{3}{c}{decomposed} \\
\cmidrule(lr){3-5}\cmidrule(lr){6-8}
$n_r$ & $N$ & $n_\mathrm{par}$~[M] & $t_\mathrm{train}$~[h] & conv. & $n_\mathrm{par}$~[M] & $t_\mathrm{train}$~[h] & conv. \\
\midrule
2 & 10 & 1.84 & 0.11 & 1.000 & 1.47 & 0.14 & 1.000 \\
2 & 20 & 2.62 & 0.23 & 0.945 & 1.47 & 0.27 & 1.000 \\
2 & 50 & 4.95 & 0.54 & 0.000 & 1.47 & 0.77 & 1.000 \\
2 & 100 & -- & -- & -- & 1.47 & 2.05 & 1.000 \\
3 & 10 & 2.23 & 0.12 & 0.996 & 1.67 & 0.15 & 1.000 \\
3 & 20 & 3.40 & 0.25 & 0.000 & 1.67 & 0.29 & 1.000 \\
3 & 50 & 6.91 & 0.61 & 0.000 & 1.67 & 0.85 & 1.000 \\
3 & 100 & -- & -- & -- & 1.67 & 2.31 & 1.000 \\
5 & 10 & 3.02 & 0.12 & 0.000 & 2.06 & 0.21 & 1.000 \\
5 & 20 & 4.97 & 0.25 & 0.000 & 2.06 & 0.39 & 1.000 \\
5 & 50 & 10.81 & 0.99 & 0.000 & 2.06 & 1.19 & 1.000 \\
5 & 100 & -- & -- & -- & 2.06 & 2.81 & 1.000 \\
\bottomrule
\end{tabular}

\end{table*}

\begin{table*}[pos=!htbp]
\centering
\caption{CPU solve times of the decomposed solver (LISCO) and IPOPT with MA27 over 256 instances per configuration, with the primal-dual dimension $n_z$ and the median, 95th percentile and maximum solve times. The speedup is the per-sample ratio $t_\mathrm{IPOPT}/t_\mathrm{LISCO}$ on identical instances, reported by its minimum, median and maximum. Both solvers converge on every instance. The monolithic solver is omitted, as it does not converge at most configurations.}
\label{tab:horizon_runtime}
\begin{tabular}{ccrrrrrrrrrr}
\toprule
 &  & \multicolumn{1}{c}{} & \multicolumn{3}{c}{$t_\mathrm{LISCO}$~[ms]} & \multicolumn{3}{c}{$t_\mathrm{IPOPT}$~[ms]} & \multicolumn{3}{c}{speedup~[$\times$]} \\
\cmidrule(lr){4-6}\cmidrule(lr){7-9}\cmidrule(lr){10-12}
$n_r$ & $N$ & $n_z$ & med & p95 & max & med & p95 & max & min & med & max \\
\midrule
2 & 10 & 380 & 0.89 & 1.29 & 1.66 & 2.65 & 2.92 & 3.39 & 1.64 & 2.99 & 4.14 \\
2 & 20 & 760 & 1.69 & 1.99 & 6.21 & 3.47 & 3.77 & 4.42 & 0.52 & 2.03 & 2.83 \\
2 & 50 & 1900 & 3.44 & 3.93 & 4.36 & 5.79 & 6.24 & 6.99 & 1.28 & 1.68 & 2.23 \\
2 & 100 & 3800 & 6.24 & 9.91 & 10.66 & 10.97 & 11.63 & 12.00 & 1.00 & 1.72 & 2.02 \\
3 & 10 & 570 & 1.08 & 1.46 & 2.05 & 3.50 & 3.81 & 4.43 & 1.65 & 3.23 & 4.86 \\
3 & 20 & 1140 & 1.72 & 2.21 & 3.42 & 5.01 & 5.52 & 5.80 & 1.47 & 2.93 & 3.81 \\
3 & 50 & 2850 & 4.79 & 5.39 & 5.89 & 8.10 & 8.72 & 9.20 & 1.37 & 1.70 & 2.00 \\
3 & 100 & 5700 & 6.64 & 7.20 & 8.32 & 14.70 & 15.68 & 16.57 & 1.71 & 2.22 & 2.55 \\
5 & 10 & 950 & 1.80 & 2.38 & 3.79 & 5.17 & 5.58 & 6.04 & 1.39 & 2.86 & 3.77 \\
5 & 20 & 1900 & 2.27 & 2.83 & 3.34 & 9.37 & 9.80 & 13.99 & 2.67 & 4.11 & 5.44 \\
5 & 50 & 4750 & 4.81 & 6.66 & 8.03 & 16.80 & 18.35 & 19.12 & 2.24 & 3.56 & 4.55 \\
5 & 100 & 9500 & 9.29 & 11.00 & 12.05 & 31.11 & 33.81 & 34.96 & 2.55 & 3.43 & 4.26 \\
\bottomrule
\end{tabular}

\end{table*}

\subsection{Smoothing-Parameter Trade-off}
\label{sec:res_eps}

\begin{figure}[pos=!htbp]
\centering
\includegraphics[width=.9\columnwidth]{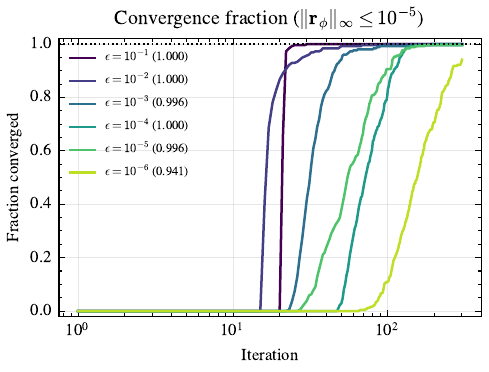}
\caption{Fraction of the 256 samples reaching the residual tolerance $\lVert \mathbf{r}_\phi \rVert_\infty \leq 10^{-5}$ over the solver iterations for the smoothing parameters $\epsilon \in \{10^{-1}, \ldots, 10^{-6}\}$ at $n_r = 2$ and $N = 50$, with an iteration cap of 300.}
\label{fig:eps_convergence_fraction}
\end{figure}

\begin{figure}[pos=!htbp]
\centering
\includegraphics[width=.9\columnwidth]{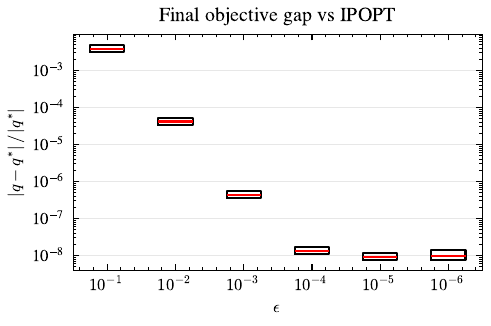}
\caption{Absolute relative objective gap $\lvert q - q^{*} \rvert / \lvert q^{*} \rvert$ at the final iterate for the smoothing sweep at $n_r = 2$ and $N = 50$ over 256 samples, where $q^{*}$ denotes the objective of the IPOPT reference solution at tolerance $10^{-8}$. The line marks the median and the boxes span the 10th to 90th percentile.}
\label{fig:eps_relative_gap}
\end{figure}

We train and evaluate the decomposed solver for the smoothing parameters $\epsilon \in \{10^{-1}, 10^{-2}, \ldots, 10^{-6}\}$ at a single problem size with $n_r = 2$ and $N = 50$, with all remaining settings fixed by Table~\ref{tab:setup}.
Fig.~\ref{fig:eps_convergence_fraction} shows the fraction of converged samples over the solver iterations for each $\epsilon$, and Fig.~\ref{fig:eps_relative_gap} shows the relative objective gap with respect to IPOPT at the final iterate.
Since IPOPT serves only as the objective reference in this study and is not included in a timing comparison, we run it at a tolerance of $10^{-8}$ without the increased barrier target of Section~\ref{sec:cs_baselines}.

The number of iterations required to reach the residual tolerance grows as $\epsilon$ decreases.
At $\epsilon = 10^{-1}$, all 256 samples converge within about 30 iterations, whereas at $\epsilon = 10^{-6}$ some runs reach the iteration cap of 300, within which 94.1\% of the samples converge.
For all intermediate values of $\epsilon$, at least 99.6\% of the samples converge. %

An important advantage of the smoothing formulation is that feasibility is not affected by the choice of $\epsilon$.
As established in Section~\ref{sec:bg_residual}, a residual below the tolerance bounds the equality constraints directly and places the inequality constraints on their feasible side for every $\epsilon$, so an accepted iterate is feasible to the tolerance level independently of the smoothing.
The smoothing relaxes the complementarity condition and thereby affects the solution optimality but not constraint satisfaction, which is an important advantage for control applications.
In Fig.~\ref{fig:eps_relative_gap}, the median relative objective gap decreases from approximately $4\cdot10^{-3}$ at $\epsilon = 10^{-1}$ to $4\cdot10^{-5}$ at $\epsilon = 10^{-2}$ and $4\cdot10^{-7}$ at $\epsilon = 10^{-3}$.
For $\epsilon \leq 10^{-4}$, the gap settles at approximately $10^{-8}$.

The smoothing parameter therefore trades optimality against convergence speed, while feasibility is retained at every $\epsilon$.
Based on this sweep, the remainder of the paper uses $\epsilon = 10^{-1}$, which yields the fastest convergence to the tolerance at a median relative objective gap of approximately $4\cdot10^{-3}$.

\subsection{Decomposition Structure and Overlap}
\label{sec:res_overlap}

\begin{figure}[pos=!htbp]
\centering
\includegraphics[width=.9\columnwidth]{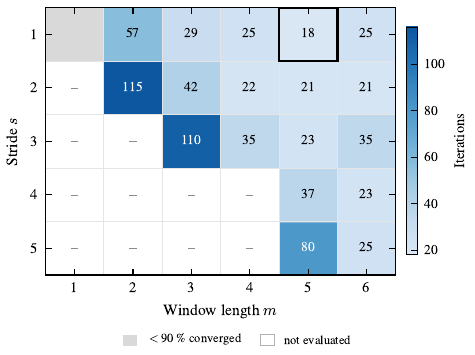}
\caption{90th percentile of the iterations to the residual tolerance $\lVert \mathbf{r}_\phi \rVert_\infty \leq 10^{-5}$ for different window lengths $m$ and strides $s$, with $n_r = 5$ reactors and $N \in \{49, 50, 51\}$. Each configuration is evaluated on 512 samples with a maximum number of 200 iterations. The grey cell indicates that fewer than 90\% of the samples converge within the iteration limit. Dashes indicate combinations that violate the admissibility condition \eqref{eq:dec_admissibility} and were not evaluated. The black outline marks the lowest observed 90th percentile.}
\label{fig:decomposition_heatmap}
\end{figure}

In the overlap study, we train and evaluate the decomposed solver for the window lengths $m \in \{1,\ldots,6\}$ and the strides $s \in \{1,\ldots,5\}$ at $n_r = 5$, $N = 50$ and $\epsilon = 10^{-1}$ and report the 90th percentile of the iterations to the residual tolerance over 512 evaluation samples in Fig.~\ref{fig:decomposition_heatmap}.
Only combinations with $s \leq m$ are admissible and the divisibility condition of \eqref{eq:dec_admissibility} must hold.
Thus, the horizon is adjusted to $N = 51$ or $N = 49$ where $N = 50$ does not satisfy it, and combinations for which no such adjustment exists are absent.

The results show that an overlap consistently lowers the iteration count.
At $m = 5$, the 90th percentile increases from 18 iterations at $s = 1$ to 80 at $s = 5$, and the non-overlapping partitions with $s = m$ are the slowest configurations that still achieve the required convergence fraction, with 115, 110 and 80 iterations for $m = 2$, $3$ and $5$.
The purely stage-local map with $m = s = 1$ does not reach 90\% convergence, as only 81.1\% of the samples converge within the cap, so window inputs beyond a single stage are necessary in this case study.
At maximum overlap, the lowest iteration count is achieved with a window length of $m = 5$.

Thus, windows of a few stages suffice for fast convergence in this case study, provided that they overlap sufficiently.

\subsection{Closed-Loop Controller Comparison}
\label{sec:res_closed_loop}

For the closed-loop comparison, we apply the decomposed solver, IPOPT and approximate MPC as NMPC controllers with $n_r = 5$ reactors and a prediction horizon of $N = 50$ steps, following Section~\ref{sec:cs_baselines}.
We simulate 50 trajectories of 50 steps each, using the same initial states for all three controllers.
These initial states are sampled uniformly over the scaling ranges in Table~\ref{tab:cstr_scaling}.
Each controller therefore executes 2500 steps, over which the mean constraint violations in Table~\ref{tbl:closed_loop_comparison} are averaged.
Constraint violations are measured on the scaled variables, so a magnitude of $10^{-1}$ corresponds to 10\% of the range of the violated variable, and Table~\ref{tbl:closed_loop_comparison} reports their mean and maximum for the input and state bounds.
Violations below $10^{-5}$ are not counted, since this is below the residual tolerance.

Table~\ref{tbl:closed_loop_comparison} lists the resulting costs, solve times and constraint violations.
The decomposed solver, trained with closed-loop resampling, reaches the residual tolerance on every executed step, and neither it nor IPOPT violates any constraint.
Its mean closed-loop cost of $-358.87$ differs from IPOPT's $-358.44$ by only 0.12\%.
The decomposed solver has a median solve time of $5.2\,\mathrm{ms}$ compared to $18.6\,\mathrm{ms}$ for IPOPT, corresponding to a speedup factor of 3.6, and its maximum of $11.6\,\mathrm{ms}$ remains below the median of IPOPT.
These times are consistent with the scaling study of Section~\ref{sec:res_scaling} at the same problem size, where the solver trained and evaluated on the uniform distribution requires a median of $4.81\,\mathrm{ms}$. %
Approximate MPC requires only $0.08\,\mathrm{ms}$ per step and achieves a similar cost to IPOPT, but produces constraint violations (Table~\ref{tbl:closed_loop_comparison}).

Approximate MPC is trained on 150000 precomputed optimizer solutions, whereas the self-supervised loss \eqref{eq:loss_full} requires no reference solutions.
Generating this training set takes 0.78 hours, which together with the 1.05 hours of supervised training gives a total offline effort of 1.83 hours for approximate MPC.
This exceeds the 0.75 hours of the complete self-supervised training of the decomposed solver.

The economic controller drives the system to a steady state on the boundary of several simultaneously active state constraints, and such states are rarely produced by uniform sampling over the scaling ranges, as discussed in Section~\ref{sec:method_cl_resampling}.
Table~\ref{tbl:step_bin_convergence} lists the fraction of solves that reach the residual tolerance within the iteration cap, binned by the simulation step, for the solvers trained with closed-loop and with uniform resampling.
On the first step, both solvers converge for all 50 uniformly sampled initial states.
From the second step onward, the trajectories leave the initial-state distribution, and the convergence of the uniformly trained solver drops to 84.5\% on steps 2 to 5, to 0.4\% on steps 6 to 10 and to zero beyond step 10, which amounts to 8.8\% over all steps.
The solver trained with closed-loop resampling converges in every bin.

\begin{table}[pos=!htbp]
\centering
\caption{Closed-loop comparison of IPOPT, the decomposed solver (LISCO) trained with closed-loop resampling and approximate MPC trained on closed-loop data, over 50 trajectories of 50 steps each, started from uniformly sampled initial states at $n_r = 5$ and $N = 50$. Reported are the mean closed-loop cost, the per-solve times and the mean and maximum constraint violations on the scaled variables for the input and state bounds, where violations below the tolerance $10^{-5}$ are not counted. The last two rows list the offline data sampling and training times.}
\label{tbl:closed_loop_comparison}
{\setlength{\tabcolsep}{3pt} \begin{tabular}{lrrr}
\toprule
metric & IPOPT & LISCO & approx. MPC \\
\midrule
closed-loop cost (mean) & -358.4415 & -358.8662 & -358.4793 \\
\midrule
solve time median [ms] & 18.61 & 5.20 & 0.08 \\
solve time p90 [ms] & 19.20 & 6.30 & 0.10 \\
solve time max [ms] & 21.41 & 11.62 & 0.43 \\
\midrule
input viol. (mean) & 0.0000 & 0.0000 & 0.0000 \\
input viol. (max) & 0.0000 & 0.0000 & 0.0002 \\
state viol. (mean) & 0.0000 & 0.0000 & 0.0004 \\
state viol. (max) & 0.0000 & 0.0000 & 0.0060 \\
\midrule
sampling time [h] & -- & -- & 0.78 \\
training time [h] & -- & 0.75 & 1.05 \\
\bottomrule
\end{tabular}
}
\end{table}

\begin{table}[pos=!htbp]
\centering
\caption{Fraction of closed-loop solves reaching the residual tolerance $\lVert \mathbf{r}_\phi \rVert_\infty \leq 10^{-5}$ within 300 iterations, binned by the simulation step, for the decomposed solver trained with closed-loop resampling ($q_{\mathrm{reset}} = 0.02$) and with uniform resampling ($q_{\mathrm{reset}} = 1$), over the 50 trajectories of 50 steps of the closed-loop study.}
\label{tbl:step_bin_convergence}
\begin{tabular}{lrr}
\toprule
step bin & closed-loop trained & uniform trained \\
\midrule
steps 1--1 & 100.0\% & 100.0\% \\
steps 2--5 & 100.0\% & 84.5\% \\
steps 6--10 & 100.0\% & 0.4\% \\
steps 11--25 & 100.0\% & 0.0\% \\
steps 26--50 & 100.0\% & 0.0\% \\
\midrule
all 1--50 & 100.0\% & 8.8\% \\
\bottomrule
\end{tabular}

\end{table}

\section{Conclusion}\label{sec:conclusion}

We introduced a learned primal-dual solver for NMPC based on temporal decomposition.
A shared neural network computes local updates in parallel on overlapping horizon segments.
This makes the network size independent of the prediction horizon while retaining self-supervised training without precomputed optimizer solutions.
We also introduced a closed-loop resampling strategy that uses predicted next states to align the training distribution with closed-loop operation at no additional sampling cost.

On an economic NMPC benchmark of a CSTR cascade with up to 100 stages and 9500 primal-dual variables, the decomposed solver converged across all tested problem sizes, including configurations where monolithic LISCO did not, and achieved median speedups of up to a factor of four over IPOPT.
With closed-loop resampling, it converged throughout the tested closed-loop trajectories and achieved control performance comparable to IPOPT with no observed constraint violations.

Future work will investigate decomposed prediction of initial iterates for warm-started operation and extend the architecture to scenario trees in robust NMPC and to weakly coupled process networks.

\printcredits

\section*{Declaration of competing interest}
The authors declare that they have no known competing financial interests or personal relationships that could have appeared to influence the work reported in this paper.

\section*{Declaration of generative AI and AI-assisted technologies in the writing process}
During the preparation of this work the authors used AI technologies (ChatGPT 5.6, Claude Fable 5 and Opus 5) to assist with grammar, phrasing, and code development. After using these tools, the authors reviewed and edited the content as needed and take full responsibility for the content of the publication.

\section*{Data availability}
The code necessary to reproduce the reported results is openly available under \url{https://github.com/lukaslueken/lisco-decompose}.

\bibliographystyle{model1-num-names}

\bibliography{references.bib}

@article{adamekEnablingRobustMixedinteger2026,
  title = {Enabling Robust Mixed-Integer Nonlinear Model Predictive Control via Self-Supervised Learning and Combinatorial Integral Approximation},
  author = {Adamek, Joshua and L{\"u}ken, Lukas and Lucia, Sergio},
  year = 2026,
  month = mar,
  journal = {J. Process Control},
  volume = {159},
  pages = {103636},
  issn = {0959-1524},
  doi = {10.1016/j.jprocont.2026.103636},
  urldate = {2026-04-29}
}

@article{anderssonCasADiSoftwareFramework2019,
  title = {{{CasADi}}: A Software Framework for Nonlinear Optimization and Optimal Control},
  shorttitle = {{{CasADi}}},
  author = {Andersson, Joel A. E. and Gillis, Joris and Horn, Greg and Rawlings, James B. and Diehl, Moritz},
  year = 2019,
  month = mar,
  journal = {Math. Program. Comput.},
  volume = {11},
  number = {1},
  pages = {1--36},
  issn = {1867-2949, 1867-2957},
  doi = {10.1007/s12532-018-0139-4},
  urldate = {2021-06-24},
  langid = {english}
}

@article{bartaliSchwarzDecompositionParallel2024,
  title = {Schwarz Decomposition for Parallel Minimum Lap-Time Problems: Evaluating against {{ADMM}}},
  author = {Bartali, Lorenzo and Gabiccini, Marco and Wright, Stephen J.},
  year = 2024,
  journal = {Veh. Syst. Dyn.},
  volume = {64},
  number = {1},
  eprint = {https://doi.org/10.1080/00423114.2024.2407175},
  pages = {1--26},
  publisher = {Taylor \& Francis},
  doi = {10.1080/00423114.2024.2407175}
}

@article{bemporadExplicitLinearQuadratic2002,
  title = {The Explicit Linear Quadratic Regulator for Constrained Systems},
  author = {Bemporad, Alberto and Morari, Manfred and Dua, Vivek and Pistikopoulos, Efstratios N.},
  year = 2002,
  month = jan,
  journal = {Automatica},
  volume = {38},
  number = {1},
  pages = {3--20},
  issn = {0005-1098},
  doi = {10.1016/S0005-1098(01)00174-1},
  urldate = {2022-11-10},
  langid = {english}
}

@article{boydDistributedOptimizationStatistical2011,
  title = {Distributed Optimization and Statistical Learning via the Alternating Direction Method of Multipliers},
  author = {Boyd, Stephen and Parikh, Neal and Chu, Eric and Peleato, Borja and Eckstein, Jonathan},
  year = 2011,
  journal = {Found. Trends Mach. Learn.},
  volume = {3},
  number = {1},
  pages = {1--122},
  publisher = {Now Publishers, Inc.}
}

@inproceedings{chenApproximatingExplicitModel2018,
  title = {Approximating Explicit Model Predictive Control Using Constrained Neural Networks},
  booktitle = {2018 {{Annual American Control Conference}} ({{ACC}})},
  author = {Chen, Steven and Saulnier, Kelsey and Atanasov, Nikolay and Lee, Daniel D. and Kumar, Vijay and Pappas, George J. and Morari, Manfred},
  year = 2018,
  month = jun,
  pages = {1520--1527},
  issn = {2378-5861},
  doi = {10.23919/ACC.2018.8431275}
}

@article{chenLearningOptimizePrimer2022,
  title = {Learning to Optimize: {{A}} Primer and a Benchmark},
  author = {Chen, Tianlong and Chen, Xiaohan and Chen, Wuyang and Wang, Zhangyang and Heaton, Howard and Liu, Jialin and Yin, Wotao},
  year = 2022,
  journal = {J. Mach. Learn. Res.},
  volume = {23},
  number = {1},
  pages = {8562--8620},
  publisher = {JMLRORG}
}

@article{chenPenalizedFischerBurmeisterNCPfunction2000,
  title = {A Penalized {{Fischer-Burmeister NCP-function}}},
  author = {Chen, Bintong and Chen, Xiaojun and Kanzow, Christian},
  year = 2000,
  month = jun,
  journal = {Math. Program.},
  volume = {88},
  number = {1},
  pages = {211--216},
  issn = {0025-5610},
  doi = {10.1007/PL00011375},
  urldate = {2024-04-04},
  copyright = {http://www.springer.com/tdm},
  langid = {english}
}

@article{diehlRealTimeIterationScheme2005,
  title = {A {{Real-Time Iteration Scheme}} for {{Nonlinear Optimization}} in {{Optimal Feedback Control}}},
  author = {Diehl, Moritz and Bock, Hans Georg and Schl{\"o}der, Johannes P.},
  year = 2005,
  month = jan,
  journal = {SIAM J. Control Optim.},
  volume = {43},
  number = {5},
  pages = {1714--1736},
  publisher = {{Society for Industrial and Applied Mathematics}},
  issn = {0363-0129},
  doi = {10.1137/S0363012902400713},
  urldate = {2022-11-12}
}

@article{dontiDC3LearningMethod2021,
  title = {{{DC3}}: {{A}} Learning Method for Optimization with Hard Constraints},
  shorttitle = {{{DC3}}},
  author = {Donti, Priya L. and Rolnick, David and Kolter, J. Zico},
  year = 2021,
  month = apr,
  journal = {arXiv:2104.12225},
  eprint = {2104.12225},
  urldate = {2021-06-18},
  archiveprefix = {arXiv}
}

@article{duff1983multifrontal,
  title = {The Multifrontal Solution of Indefinite Sparse Symmetric Linear Equations},
  author = {Duff, Iain S and Reid, John K},
  year = 1983,
  journal = {ACM Trans. Math. Softw.},
  volume = {9},
  number = {3},
  pages = {302--325},
  publisher = {ACM New York, NY, USA}
}

@article{ellisTutorialReviewEconomic2014,
  title = {A Tutorial Review of Economic Model Predictive Control Methods},
  author = {Ellis, Matthew and Durand, Helen and Christofides, Panagiotis D.},
  year = 2014,
  month = aug,
  journal = {J. Process Control},
  series = {Economic Nonlinear Model Predictive Control},
  volume = {24},
  number = {8},
  pages = {1156--1178},
  issn = {0959-1524},
  doi = {10.1016/j.jprocont.2014.03.010},
  urldate = {2026-08-12}
}

@article{ferreauEmbeddedOptimizationMethods2017,
  title = {Embedded {{Optimization Methods}} for {{Industrial Automatic Control}}},
  author = {Ferreau, H. J. and Alm{\'e}r, S. and Verschueren, R. and Diehl, M. and Frick, D. and Domahidi, A. and Jerez, J. L. and Stathopoulos, G. and Jones, C.},
  year = 2017,
  month = jul,
  journal = {IFAC-PapersOnLine},
  series = {20th {{IFAC World Congress}}},
  volume = {50},
  number = {1},
  pages = {13194--13209},
  issn = {2405-8963},
  doi = {10.1016/j.ifacol.2017.08.1946},
  urldate = {2022-12-21},
  langid = {english}
}

@inproceedings{fiorettoLagrangianDualityConstrained2020,
  title = {Lagrangian Duality for Constrained Deep Learning},
  booktitle = {Joint {{European Conference}} on {{Machine Learning}} and {{Knowledge Discovery}} in {{Databases}}},
  author = {Fioretto, Ferdinando and Hentenryck, Pascal Van and Mak, Terrence WK and Tran, Cuong and Baldo, Federico and Lombardi, Michele},
  year = 2020,
  pages = {118--135},
  publisher = {Springer}
}

@article{fischerSpecialNewtontypeOptimization1992,
  title = {A Special Newton-Type Optimization Method},
  author = {Fischer, A.},
  year = 1992,
  month = jan,
  journal = {Optimization},
  volume = {24},
  number = {3-4},
  pages = {269--284},
  publisher = {Taylor \& Francis},
  issn = {0233-1934},
  doi = {10.1080/02331939208843795},
  urldate = {2023-10-13}
}

@article{frisonHPIPMHighperformanceQuadratic2020,
  title = {{{HPIPM}}: A High-Performance Quadratic Programming Framework for Model Predictive Control},
  shorttitle = {{{HPIPM}}},
  author = {Frison, Gianluca and Diehl, Moritz},
  year = 2020,
  month = jan,
  journal = {IFAC-PapersOnLine},
  series = {21st {{IFAC World Congress}}},
  volume = {53},
  number = {2},
  pages = {6563--6569},
  issn = {2405-8963},
  doi = {10.1016/j.ifacol.2020.12.073},
  urldate = {2023-07-26},
  langid = {english}
}

@inproceedings{frommeUnsupervisedClosedloopPrimaldual2025,
  title = {Unsupervised Closed-Loop Primal-Dual Learning of Approximate Model Predictive Controllers},
  booktitle = {2025 European Control Conference ({{ECC}})},
  author = {Fromme, Florian and L{\"u}ken, Lukas and Lucia, Sergio},
  year = 2025,
  pages = {2970--2976},
  publisher = {IEEE}
}

@article{heinleinRobustModelPredictive2025,
  title = {Robust {{Model Predictive Control Exploiting Monotonicity Properties}}},
  author = {Heinlein, Moritz and Subramanian, Sankaranarayanan and Lucia, Sergio},
  year = 2025,
  month = sep,
  journal = {IEEE Trans. Autom. Control},
  volume = {70},
  number = {9},
  pages = {6260--6267},
  issn = {1558-2523},
  doi = {10.1109/TAC.2025.3558137},
  urldate = {2025-12-09}
}

@article{hendrycksGaussianErrorLinear2016,
  title = {Gaussian Error Linear Units (Gelus)},
  author = {Hendrycks, Dan and Gimpel, Kevin},
  year = 2016,
  journal = {arXiv:1606.08415},
  eprint = {1606.08415},
  archiveprefix = {arXiv}
}

@article{hertneckLearningApproximateModel2018,
  title = {Learning an Approximate Model Predictive Controller with Guarantees},
  author = {Hertneck, Michael and K{\"o}hler, Johannes and Trimpe, Sebastian and Allg{\"o}wer, Frank},
  year = 2018,
  month = jul,
  journal = {IEEE Control Syst. Lett.},
  volume = {2},
  number = {3},
  pages = {543--548},
  issn = {2475-1456},
  doi = {10.1109/LCSYS.2018.2843682}
}

@misc{jeonCusADiGPUParallelization2024,
  title = {{{CusADi}}: {{A GPU Parallelization Framework}} for {{Symbolic Expressions}} and {{Optimal Control}}},
  shorttitle = {{{CusADi}}},
  author = {Jeon, Se Hwan and Hong, Seungwoo and Lee, Ho Jae and Khazoom, Charles and Kim, Sangbae},
  year = 2024,
  month = aug,
  journal = {arXiv.org},
  urldate = {2025-11-24},
  howpublished = {https://arxiv.org/abs/2408.09662v1},
  langid = {english}
}

@article{kargEfficientRepresentationApproximation2020,
  title = {Efficient Representation and Approximation of Model Predictive Control Laws via Deep Learning},
  author = {Karg, Benjamin and Lucia, Sergio},
  year = 2020,
  month = sep,
  journal = {IEEE Trans. Cybern.},
  volume = {50},
  number = {9},
  pages = {3866--3878},
  issn = {2168-2267, 2168-2275},
  doi = {10.1109/TCYB.2020.2999556},
  urldate = {2021-02-16},
  langid = {english}
}

@article{kargProbabilisticPerformanceValidation2021,
  title = {Probabilistic Performance Validation of Deep Learning-Based Robust {{NMPC}} Controllers},
  author = {Karg, Benjamin and Alamo, Teodoro and Lucia, Sergio},
  year = 2021,
  journal = {Int. J. Robust Nonlinear Control},
  volume = {31},
  number = {18},
  pages = {8855--8876},
  issn = {1099-1239},
  doi = {10.1002/rnc.5696},
  urldate = {2022-11-18},
  langid = {english}
}

@article{klaucoMachineLearningbasedWarm2019,
  title = {Machine Learning-Based Warm Starting of Active Set Methods in Embedded Model Predictive Control},
  author = {Klau{\v c}o, Martin and Kal{\'u}z, Martin and Kvasnica, Michal},
  year = 2019,
  month = jan,
  journal = {Eng. Appl. Artif. Intell.},
  volume = {77},
  pages = {1--8},
  issn = {0952-1976},
  doi = {10.1016/j.engappai.2018.09.014},
  urldate = {2023-09-11}
}

@article{kotaryEndtoEndConstrainedOptimization2021,
  title = {End-to-{{End Constrained Optimization Learning}}: {{A Survey}}},
  shorttitle = {End-to-{{End Constrained Optimization Learning}}},
  author = {Kotary, James and Fioretto, Ferdinando and Van Hentenryck, Pascal and Wilder, Bryan},
  year = 2021,
  month = mar,
  journal = {arXiv:2103.16378 [cs]},
  eprint = {2103.16378},
  primaryclass = {cs},
  urldate = {2022-04-21},
  archiveprefix = {arXiv}
}

@inproceedings{lamNumbaLLVMbasedPython2015,
  title = {Numba: A {{LLVM-based Python JIT}} Compiler},
  shorttitle = {Numba},
  booktitle = {Proceedings of the {{Second Workshop}} on the {{LLVM Compiler Infrastructure}} in {{HPC}}},
  author = {Lam, Siu Kwan and Pitrou, Antoine and Seibert, Stanley},
  year = 2015,
  month = nov,
  series = {{{LLVM}} '15},
  pages = {1--6},
  publisher = {Association for Computing Machinery},
  address = {New York, NY, USA},
  doi = {10.1145/2833157.2833162},
  urldate = {2026-08-13},
  isbn = {978-1-4503-4005-2}
}

@article{liao-mcphersonFBstabProximallyStabilized2020,
  title = {{{FBstab}}: {{A}} Proximally Stabilized Semismooth Algorithm for Convex Quadratic Programming},
  shorttitle = {{{FBstab}}},
  author = {{Liao-McPherson}, Dominic and Kolmanovsky, Ilya},
  year = 2020,
  month = mar,
  journal = {Automatica},
  volume = {113},
  pages = {108801},
  issn = {0005-1098},
  doi = {10.1016/j.automatica.2019.108801},
  urldate = {2023-12-18}
}

@article{liao-mcphersonRegularizedSmoothedFischer2019,
  title = {A {{Regularized}} and {{Smoothed Fischer}}--{{Burmeister Method}} for {{Quadratic Programming With Applications}} to {{Model Predictive Control}}},
  author = {{Liao-McPherson}, Dominic and Huang, Mike and Kolmanovsky, Ilya},
  year = 2019,
  month = jul,
  journal = {IEEE Trans. Autom. Control},
  volume = {64},
  number = {7},
  pages = {2937--2944},
  issn = {1558-2523},
  doi = {10.1109/TAC.2018.2872201}
}

@article{liLearningOptimize2016,
  title = {Learning to Optimize},
  author = {Li, Ke and Malik, Jitendra},
  year = 2016,
  month = jun,
  journal = {arXiv:1606.01885},
  eprint = {1606.01885},
  urldate = {2022-03-21},
  archiveprefix = {arXiv},
  langid = {english}
}

@inproceedings{liPDHGunrolledLearningtooptimizeMethod2024,
  title = {{{PDHG-unrolled}} Learning-to-Optimize Method for Large-Scale Linear Programming},
  booktitle = {Proceedings of the 41st International Conference on Machine Learning},
  author = {Li, Bingheng and Yang, Linxin and Chen, Yupeng and Wang, Senmiao and Mao, Haitao and Chen, Qian and Ma, Yao and Wang, Akang and Ding, Tian and Tang, Jiliang and Sun, Ruoyu},
  editor = {Salakhutdinov, Ruslan and Kolter, Zico and Heller, Katherine and Weller, Adrian and Oliver, Nuria and Scarlett, Jonathan and Berkenkamp, Felix},
  year = 2024,
  month = jul,
  series = {Proceedings of Machine Learning Research},
  volume = {235},
  pages = {29164--29180},
  publisher = {PMLR}
}

@article{loshchilov2017decoupled,
  title = {Decoupled Weight Decay Regularization},
  author = {Loshchilov, Ilya and Hutter, Frank},
  year = 2017,
  journal = {arXiv preprint arXiv:1711.05101},
  eprint = {1711.05101},
  archiveprefix = {arXiv}
}

@article{luciaMultistageNonlinearModel2013,
  title = {Multi-Stage Nonlinear Model Predictive Control Applied to a Semi-Batch Polymerization Reactor under Uncertainty},
  author = {Lucia, Sergio and Finkler, Tiago and Engell, Sebastian},
  year = 2013,
  journal = {J. Process Control},
  volume = {23},
  number = {9},
  pages = {1306--1319},
  issn = {0959-1524},
  doi = {10.1016/j.jprocont.2013.08.008},
  urldate = {2022-02-10},
  langid = {english}
}

@inproceedings{lukenLearningIterativeSolvers2024,
  title = {Learning Iterative Solvers for Accurate and Fast Nonlinear Model Predictive Control via Unsupervised Training},
  booktitle = {2024 European Control Conference ({{ECC}})},
  author = {L{\"u}ken, Lukas and Lucia, Sergio},
  year = 2024,
  pages = {1843--1850},
  doi = {10.23919/ECC64448.2024.10591099}
}

@article{lukenSelfsupervisedLearningIterative2026,
  title = {Self-Supervised Learning of Iterative Solvers for Constrained Optimization},
  author = {L{\"u}ken, Lukas and Lucia, Sergio},
  year = 2026,
  month = jun,
  journal = {Results Control Optim.},
  volume = {23},
  pages = {100751},
  issn = {2666-7207},
  doi = {10.1016/j.rico.2026.100751},
  urldate = {2026-07-01}
}

@article{lukenSobolevTrainingDataefficient2023,
  title = {Sobolev Training for Data-Efficient Approximate Nonlinear {{MPC}}},
  author = {L{\"u}ken, Lukas and Brandner, Dean and Lucia, Sergio},
  year = 2023,
  journal = {IFAC-PapersOnLine},
  volume = {56},
  number = {2},
  pages = {5765--5772},
  publisher = {Elsevier}
}

@article{mitsosOptimalDeterministicAlgorithm2018,
  title = {Optimal Deterministic Algorithm Generation},
  author = {Mitsos, Alexander and Najman, Jaromi{\l} and Kevrekidis, Ioannis G.},
  year = 2018,
  month = aug,
  journal = {J. Glob. Optim.},
  volume = {71},
  number = {4},
  pages = {891--913},
  issn = {1573-2916},
  doi = {10.1007/s10898-018-0611-8},
  urldate = {2023-05-15},
  langid = {english}
}

@article{naSuperconvergenceOnlineOptimization2023,
  title = {Superconvergence of {{Online Optimization}} for {{Model Predictive Control}}},
  author = {Na, Sen and Anitescu, Mihai},
  year = 2023,
  month = mar,
  journal = {IEEE Trans. Autom. Control},
  volume = {68},
  number = {3},
  pages = {1383--1398},
  issn = {1558-2523},
  doi = {10.1109/TAC.2022.3223323},
  urldate = {2024-08-29}
}

@book{nocedalNumericalOptimization2006,
  title = {Numerical Optimization},
  author = {Nocedal, Jorge and Wright, Stephen J.},
  year = 2006,
  series = {Springer Series in Operations Research},
  edition = {2nd},
  publisher = {Springer},
  address = {New York},
  isbn = {978-0-387-30303-1},
  langid = {english},
  lccn = {QA402.5 .N62 2006}
}

@article{parisiniRecedinghorizonRegulatorNonlinear1995,
  title = {A Receding-Horizon Regulator for Nonlinear Systems and a Neural Approximation},
  author = {Parisini, T. and Zoppoli, R.},
  year = 1995,
  month = oct,
  journal = {Automatica},
  volume = {31},
  number = {10},
  pages = {1443--1451},
  issn = {0005-1098},
  doi = {10.1016/0005-1098(95)00044-W},
  urldate = {2022-11-12},
  langid = {english}
}

@article{parkSelfsupervisedPrimaldualLearning2023,
  title = {Self-Supervised Primal-Dual Learning for Constrained Optimization},
  author = {Park, Seonho and Van Hentenryck, Pascal},
  year = 2023,
  month = jun,
  journal = {Proc. AAAI Conf. Artif. Intell.},
  volume = {37},
  number = {4},
  pages = {4052--4060},
  issn = {2374-3468, 2159-5399},
  doi = {10.1609/aaai.v37i4.25520},
  urldate = {2024-01-10},
  langid = {english}
}

@inproceedings{paszkePyTorchImperativeStyle2019,
  title = {{{PyTorch}}: {{An}} Imperative Style, High-Performance Deep Learning Library},
  booktitle = {Advances in Neural Information Processing Systems},
  author = {Paszke, Adam and Gross, Sam and Massa, Francisco and Lerer, Adam and Bradbury, James and Chanan, Gregory and Killeen, Trevor and Lin, Zeming and Gimelshein, Natalia and Antiga, Luca and Desmaison, Alban and Kopf, Andreas and Yang, Edward and DeVito, Zachary and Raison, Martin and Tejani, Alykhan and Chilamkurthy, Sasank and Steiner, Benoit and Fang, Lu and Bai, Junjie and Chintala, Soumith},
  editor = {Wallach, H. and Larochelle, H. and Beygelzimer, A. and {dAlch{\'e}-Buc}, F. and Fox, E. and Garnett, R.},
  year = 2019,
  volume = {32},
  publisher = {Curran Associates, Inc.}
}

@inproceedings{rawlingsFundamentalsEconomicModel2012a,
  title = {Fundamentals of Economic Model Predictive Control},
  booktitle = {2012 {{IEEE}} 51st {{IEEE Conference}} on {{Decision}} and {{Control}} ({{CDC}})},
  author = {Rawlings, James B. and Angeli, David and Bates, Cuyler N.},
  year = 2012,
  month = dec,
  pages = {3851--3861},
  issn = {0743-1546},
  doi = {10.1109/CDC.2012.6425822},
  urldate = {2026-08-12}
}

@book{rawlingsModelPredictiveControl2017,
  title = {Model Predictive Control: Theory, Computation, and Design},
  shorttitle = {Model Predictive Control},
  author = {Rawlings, James B. and Mayne, David Q. and Diehl, Moritz M.},
  year = 2017,
  edition = {2nd},
  publisher = {Nob Hill Publishing},
  address = {Madison, Wisconsin},
  isbn = {978-0-9759377-3-0},
  langid = {english}
}

@article{sambharya2026learning,
  title = {Learning Algorithm Hyperparameters for Fast Parametric Convex Optimization},
  author = {Sambharya, Rajiv and Stellato, Bartolomeo},
  year = 2026,
  journal = {SIAM J. Math. Data Sci.},
  volume = {8},
  number = {3},
  pages = {649--676},
  publisher = {SIAM}
}

@inproceedings{sambharyaEndtoendLearningWarmstart2023,
  title = {End-to-End Learning to Warm-Start for Real-Time Quadratic Optimization},
  booktitle = {Learning for Dynamics and Control Conference},
  author = {Sambharya, Rajiv and Hall, Georgina and Amos, Brandon and Stellato, Bartolomeo},
  year = 2023,
  pages = {220--234},
  publisher = {PMLR}
}

@article{shinExponentialDecaySensitivity2022,
  title = {Exponential {{Decay}} of {{Sensitivity}} in {{Graph-Structured Nonlinear Programs}}},
  author = {Shin, Sungho and Anitescu, Mihai and Zavala, Victor M.},
  year = 2022,
  month = jun,
  journal = {SIAM J. Optim.},
  volume = {32},
  pages = {1156--1183},
  publisher = {{Society for Industrial and Applied Mathematics}},
  issn = {1052-6234},
  doi = {10.1137/21M1391079},
  urldate = {2022-06-09}
}

@article{shlezingerModelBasedDeepLearning2023,
  title = {Model-{{Based Deep Learning}}},
  author = {Shlezinger, Nir and Whang, Jay and Eldar, Yonina C. and Dimakis, Alexandros G.},
  year = 2023,
  month = may,
  journal = {Proc. IEEE},
  volume = {111},
  number = {5},
  pages = {465--499},
  issn = {1558-2256},
  doi = {10.1109/JPROC.2023.3247480}
}

@article{vaupelAcceleratingNonlinearModel2020,
  title = {Accelerating Nonlinear Model Predictive Control through Machine Learning},
  author = {Vaupel, Yannic and Hamacher, Nils C. and Caspari, Adrian and Mhamdi, Adel and Kevrekidis, Ioannis G. and Mitsos, Alexander},
  year = 2020,
  month = aug,
  journal = {J. Process Control},
  volume = {92},
  pages = {261--270},
  issn = {09591524},
  doi = {10.1016/j.jprocont.2020.06.012},
  urldate = {2023-09-13},
  langid = {english}
}

@article{verschuerenAcadosModularOpensource2022,
  title = {Acados---a Modular Open-Source Framework for Fast Embedded Optimal Control},
  author = {Verschueren, Robin and Frison, Gianluca and Kouzoupis, Dimitris and Frey, Jonathan and van Duijkeren, Niels and Zanelli, Andrea and Novoselnik, Branimir and Albin, Thivaharan and Quirynen, Rien and Diehl, Moritz},
  year = 2022,
  month = mar,
  journal = {Math. Program. Comput.},
  volume = {14},
  number = {1},
  pages = {147--183},
  issn = {1867-2957},
  doi = {10.1007/s12532-021-00208-8},
  urldate = {2023-11-21},
  langid = {english}
}

@article{wachterImplementationInteriorpointFilter2006,
  title = {On the Implementation of an Interior-Point Filter Line-Search Algorithm for Large-Scale Nonlinear Programming},
  author = {W{\"a}chter, Andreas and Biegler, Lorenz T.},
  year = 2006,
  month = mar,
  journal = {Math. Program.},
  volume = {106},
  number = {1},
  pages = {25--57},
  issn = {1436-4646},
  doi = {10.1007/s10107-004-0559-y},
  urldate = {2022-04-25},
  langid = {english}
}

@article{zanelli_forces_2020,
  title = {{{FORCES NLP}}: An Efficient Implementation of Interior-Point Methods for Multistage Nonlinear Nonconvex Programs},
  shorttitle = {{{FORCES NLP}}},
  author = {Zanelli, Andrea and Domahidi, Alexander and Jerez, Juan and Morari, Manfred},
  year = 2020,
  journal = {Int. J. Control},
  volume = {93},
  number = {1},
  pages = {13--29}
}

\end{document}